\documentclass[12pt,a4paper,twoside]{article}
\usepackage[T1]{fontenc}
\usepackage{amsfonts} 
\usepackage{amssymb}
\usepackage{amsmath}
\usepackage{url}
\usepackage{latexsym}

\usepackage{mathbbol} 
\usepackage{epsfig}
\usepackage{placeins}  
\usepackage{mathrsfs} 

\newcommand{\1}{\Eins}
\newcommand{\sqc}{\mathrel{\mathrm{X}\mkern-14mu\square}}
\newcommand{\Depth}{\operatorname{Depth}}
\newcommand{\Height}{\operatorname{Height}}
\newcommand{\Sol}{\mathfrak{S}}
\newcommand{\F}{\mathscr{F}}
\newcommand{\Fcal}{\mathcal{F}}
\newcommand{\pfrak}{\mathfrak{p}}
\newcommand{\D}{\mathscr{D}}

\newcommand{\cg}{{\mathfrak{f}}} 
\newcommand{\corr}{\mathscr{C}}
\newcommand{\cvv}[1]{\check{#1}}
\newcommand{\voisin}{\mathrel{\hbox to 0pt{$\mkern4mu\circ$\hss}%
\raise2.3pt\hbox{\vrule height0.2pt depth0.1pt width 9pt}}}

\newtheorem{thm}{Theorem}[section]
\newtheorem{lem}{Lemma}[section]

\newtheorem{prop}{Proposition}[section]
\newtheorem{conj}{Conjecture}[section]
\newtheorem{defi}{Definition}[section]
\newenvironment{dem}{\medbreak\noindent{\sl Proof: }}%
{\hfill $\diamond\diamond\diamond$\par\medbreak}

\begin{document}

\title{A stronger model for peg solitaire, II\footnote{keywords: Peg
    solitaire, Hi-Q, Pagoda function.}{}\footnote{AMS classification:
    primary 05A99, secondary  91A46, 52B12, 90C08.}}
\author{O. Ramar\'e,\\ CNRS, Laboratoire Painlev\'e, Universit\'e
  Lille 1\\ 59\ 655 Villeneuve d'Ascq, France}

\maketitle

\begin{abstract}
  The main problem addressed here is to decide whether it is or not possible
  to go from a given position on a peg-solitaire board to another one. No
  non-trivial sufficient conditions are known, but tests have been devised to
  show it is not possible. We expose the way these tests work in a unified
  formalism and provide a new one which is strictly stronger than all previous ones.
\end{abstract}

\section{Introduction}

Peg solitaire (also called Hi-Q) is a very simple board game that
appeared in Europe most probably at the end of the 17th century. Its
prior origin is unknown. The first evidence is a painting by
Claude-Auguste Berey of Anne Chabot de Rohan (1663-1709) playing it.
It seems to have then become popular in some royal courts. The
mathematical study of the game starts in 1710 when Leibniz writes a
memoir on the subject~\cite{Leibniz*10}.  We refer the reader to the
excellent historical account presented in Beasley's
book~\cite{Beasley*92}.  Let us introduce rapidly how this game is
being played. The first data is a board $\Sol$ which in first
approximation may be thought of as a subset of $\mathbb{Z}^2$. The
classical ones are the english board and the french one drawn below,
and we present a third one introduced by J.C. Wiegleb in~1779 (see
\cite{Beasley*92}).

\begin{figure}[!h]
  \centering
  \begin{minipage}[c]{0.26\textwidth}
     \centering \scalebox{0.20}{\includegraphics{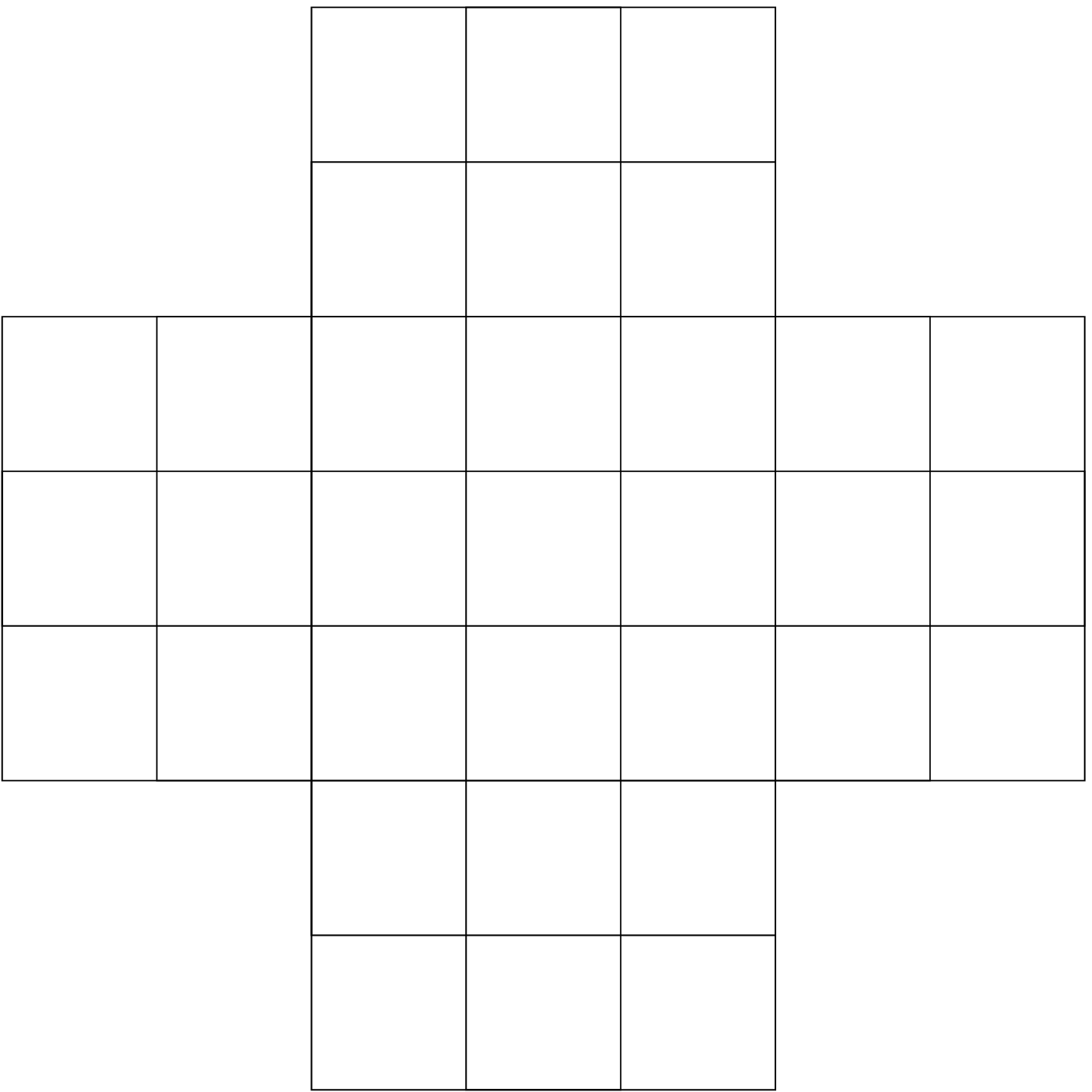}}%
     \caption{English board}
  \end{minipage}\quad\quad
  \begin{minipage}[c]{0.26\textwidth}
     \centering \scalebox{0.20}{\includegraphics{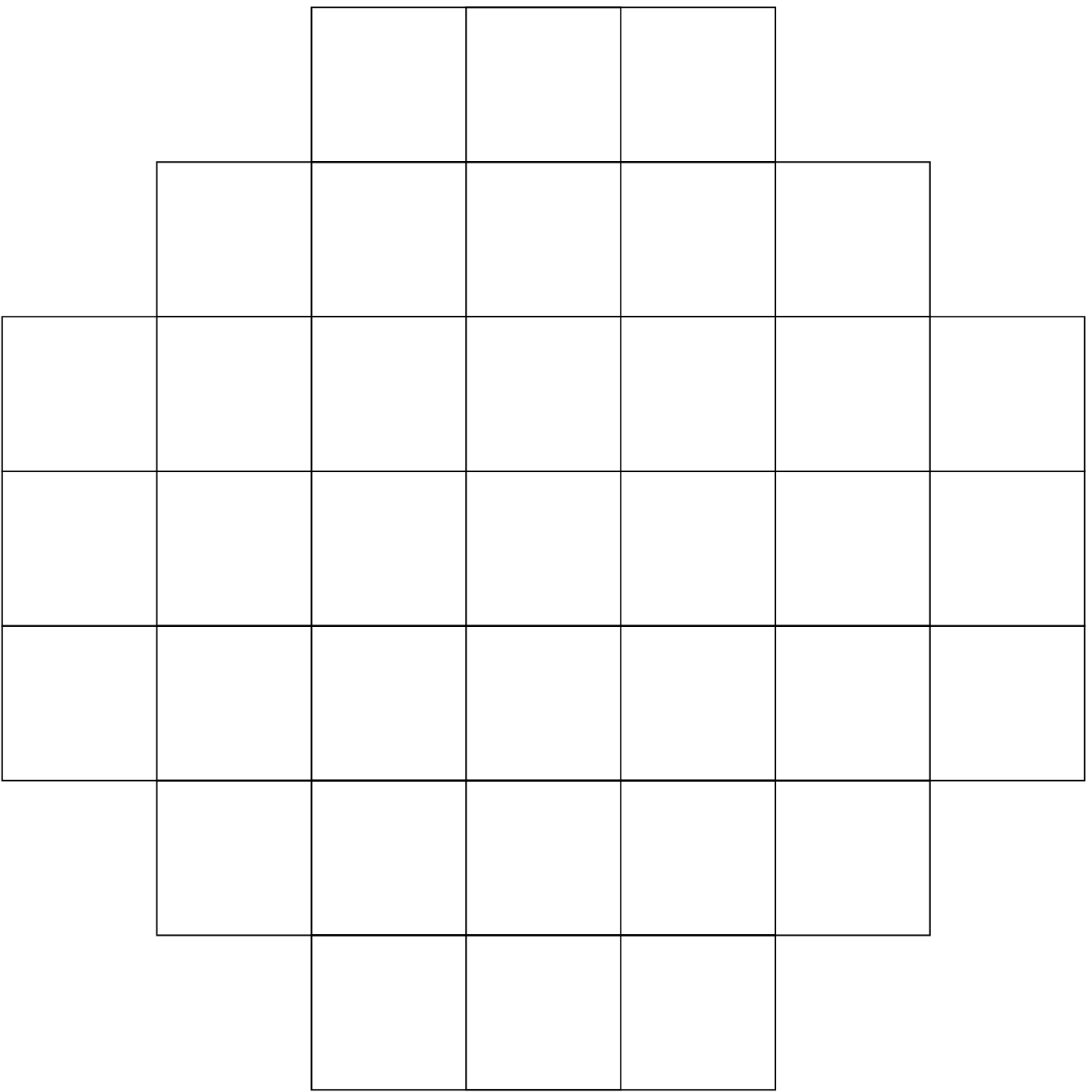}}%
     \caption{French board}
  \end{minipage}\quad\quad
  \begin{minipage}[c]{0.26\textwidth}
     \centering \scalebox{0.20}{\includegraphics{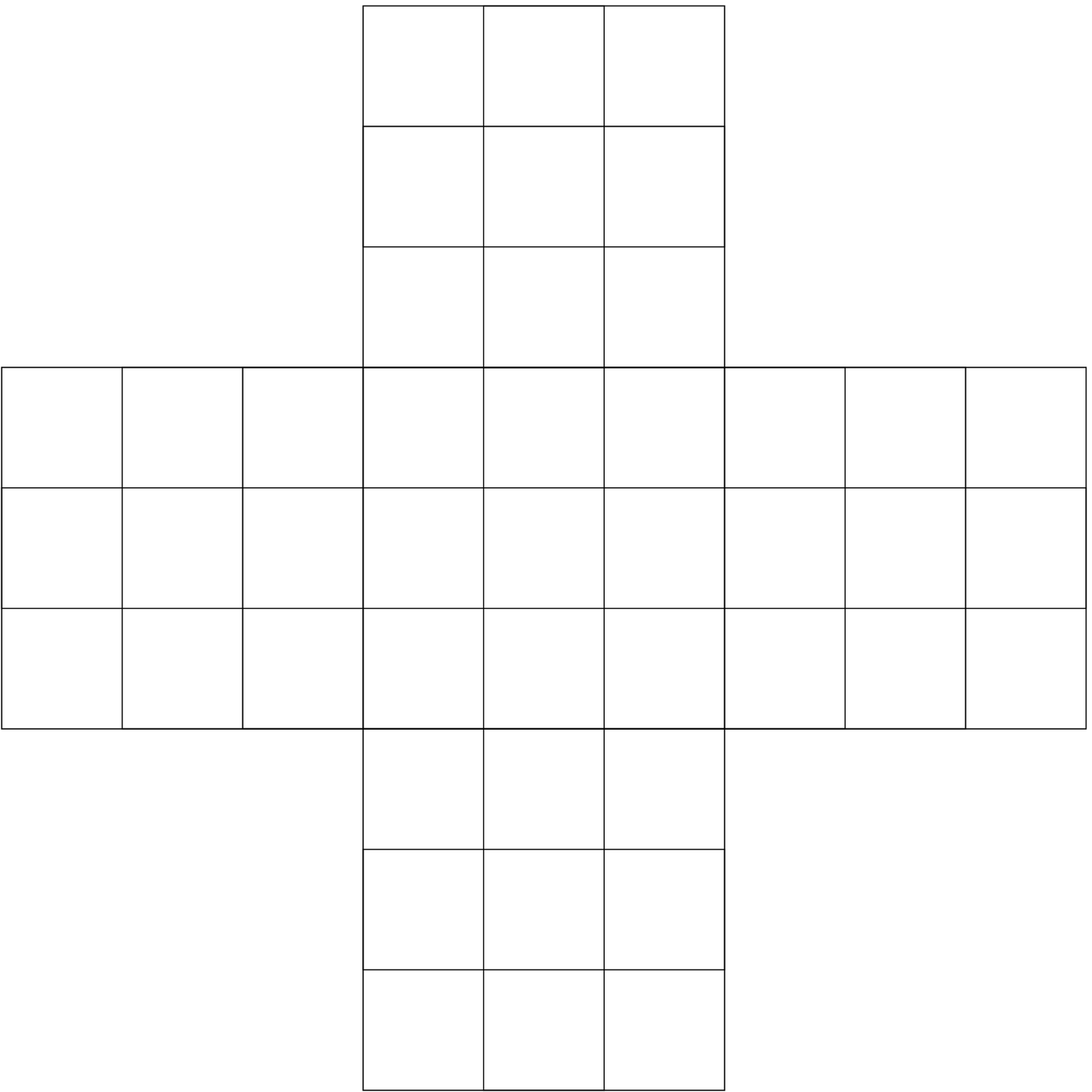}}%
     \caption{Wiegleb board}
  \end{minipage}
\end{figure}


Each square of this board can hold at most one peg, and a \emph{problem} as
we define it here is to go from a given distribution of these pegs (say $I$)
to another one (say $J$), via a succession of \emph{legal moves} that we now
define.  Given three consecutive squares $P$, $Q$ and $R$ in a row or a column
(but \emph{not} on a diagonal), of which two consecutive (say $P$ and $Q$)
contain a peg while the third one ($R$) does not, a \emph{legal move} consists
in removing the two pegs in $P$ and $Q$ and putting one on the empty square
$R$. We classically say that the closest peg in $P$ jumps over the middle one
in $Q$ and lands in $R$, while destroying the peg in $Q$. As a trivial
consequence, the number of pegs on the board decreases when the game proceeds
further. For most authors, a problem consists in reducing the initial
distribution of pegs, what we call thereafter the \emph{initial position}, to
a single peg via legal moves. They qualify the position as \emph{solvable} if
this is possible. We shall say that the problem in our sense is
\emph{feasible} if one can go from the initial position to the final one by
using legal moves. Note that the number of such moves is known and equals the
difference between the number of pegs in the initial position and the number
of pegs in the final one (that is: $|I|-|J|$).

Given a problem, we can try all possible legal moves and repeat this action until
the required number ($|I|-|J|$) of moves is reached or no further move is
possible. This process usually gets stuck because of the combinatorial
explosion. For instance E. Harang \cite{Harang*97}
computed that there are $577\ 116\ 156\ 815\ 309\ 849\ 672$ paths
on the english board from the initial position consisting of the full board on
which we leave the central square empty. Of which 
$40\ 861\ 647\ 040\ 079\ 968$ lead to the final peg being on the central square.
See also \cite{DosSantos*99}. 
In fact, numerous setting tend to show that the problem is NP-complete. For
this sentence to have a sense, we are to choose a way of extending the board to
infinity, and there is no canonical fashion to achieve that. The case of an
$n\times n$ board is studied in \cite{Uehara-Iwata*90} while the $k\times n$
board with $k$ fixed is shown to be linear in \cite{Ravikumar*97}. Of course,
one may wonder whether the english board as a subset of a $7\times7$ board is
tractable or not and the answer is still no, at least not without huge resources.
The number of paths being enormous, we look for tests
that will ensure us that it is not possible to solve a given problem. We would
welcome any test that would guarantee the feasibility, but none are yet known.

The first of this test is attributed to Reiss in~1857 in \cite{Reiss*57}
though Beasley traces it back to A. Suremain de Missery, a former officer of
the French artillery, around~1842. We again refer the reader
to~\cite{Beasley*92} for more historical details. It is also described in
Lucas book \cite{Lucas*91}, which contains also more material and in 
the dedicated chapter of~\cite{Berlekamp-Conway-Guy*82}.
A seemingly more algebraic approach is proposed in~\cite{deBruijn*72}, but it 
turns out to be only a different setting for the same test. 
This test is very often reduced by modern authors to the rule-of-three
test (see below).

We shall first present these tests in a formalism that will help us clarify the
situation; this formalism will also be adequate to present the
advances realised on the subject  in 1961/1962 at Cambridge university by
a group of students (among which were Beasley) led by J.H.~Conway.

We shall finally present a different test, which we term quadratic, and which is stronger
than all previous ones. It however relies on solving a larger integer linear
program and can sometime be resource demanding. We provide however
numerous examples that we have discovered by
exploring thousands of problems, and this in itself shows the
practicality of the approach. The theory of this test in
its purest form is complete, but we provide in the two last sections several
improvements of it, on which we are still working.
All examples have been computed via an intensive use of the lp\_solve
library~\cite{Lp-solve*06}, a GTK interface and a C-program both due to the
author. 

Let us end this introduction by mentioning that
Beasley also introduced a very geometrical tool (the \emph{in and out Theorems}),
but it does not fit well in our framework and has not been worked out 
for an arbitrary problem (to the best of my knowledge at least), even if one
remains on an english board. We shall not discuss it here.
In more recent time, there has been attempts at working out a model of this game via
string rewriting as in~\cite{Ravikumar*97}. This approach remains however
fundamentally one dimensional as are string rewriting rules. It has
had applications though in describing the complexity of the game.

\section{Main formalism of the linear board}

Given a board $\Sol$, we consider the $\mathbb{Z}$-module
$\F(\Sol,\mathbb{Z})$ of all rational integer valued functions over this
board, and define similarly $\F(\Sol,\mathbb{F}_2)$ and
$\F(\Sol,\mathbb{Q})$. This is one of the main step of the formalization: a
position in the game is given by a subset $I\subset\Sol$ (the set of squares
containing a peg), which we model by its characteristic function $\1_I$. 
If $P\in\Sol$, we note $\cvv{P}$ the function that is 1 in $P$ and 0
everywhere else. A move is thus the function $\cg=\cvv{P}+\cvv{Q}-\cvv{R}$
and $\1_I-\cg$ should become another characteristic function; we have of
course assumed that $P$, $Q$ and $R$ where three consecutive points in this
order either in a row or in a column of $\Sol$. We denote
the set of these moves by $\D(\Sol)$. In the case of the english board,
$\D(\Sol)$ has cardinality~76, while $\Sol$ has cardinality~33.

Here comes the main remark. Assume we can go from $I$ to $J$ by the succession
of legal moves $\cg_1,\cg_2,\dots, \cg_k$. Then we have
\begin{equation}
  \label{eq:1}
  \1_I-\1_J = \sum_{1\le i\le k}\cg_i.
\end{equation}
There are three ways to exploit this writing. We can say that
\begin{itemize}
\item $\1_I-\1_J$ is a rational integer linear combination of members of
  $\D(\Sol)$. This leads to the classical Reiss's theory, or to the lattice
  criterion of~\cite{Avis-Deza-Onn*00}.
\item $\1_I-\1_J$ is a linear combination with non-negative rational
  coefficients of members of  $\D(\Sol)$. This leads to the main part of
  Conway's group theory.
\item $\1_I-\1_J$ is a linear combination with non-negative integer
  coefficients of members of  $\D(\Sol)$. This leads to what we call the
  \emph{full linear test}, or also the non-negative integer test.
\end{itemize}
We introduce some notations
\begin{equation}
  \label{eq:2}
  V(\Sol,\mathbb{Z})=\sum_{\cg\in\D(\Sol)} \mathbb{Z}\cdot \cg
\end{equation}
and
\begin{equation}
  \label{eq:3}
  V^+(\Sol,\mathbb{Q})=\sum_{\cg\in\D(\Sol)} \mathbb{Q}^+\cdot \cg
  \quad,\quad
  V^+(\Sol,\mathbb{Z})=\sum_{\cg\in\D(\Sol)} \mathbb{Z}^+\cdot \cg.
\end{equation}

\section{Reiss theory and the rule-of-three test}

Let us first expose rapidly and in modern notations the classical
material. Characteristic functions having values~0 or~1, it is
tempting to look at $\1_I$ as taking its values in the field with two
elements $\mathbb {F}_2$. To avoid confusion, we note $\tilde{\1}_I$
this characteristic function as an element of $\F(\Sol,\mathbb{F}_2)$.
If one can go from the initial position $I$ to the final one $J$ by
the succession of legal moves $\cg_1,\cg_2,\dots, \cg_k$, one still
has
\begin{equation*}
   \tilde{\1}_I-\tilde{\1}_J = \sum_{1\le i\le k}\tilde{\cg}_i
\end{equation*}
where $\tilde{\cg}_i$ are of course the moves seen with values in
$\mathbb{F}_2$. If $\cg=\cvv{P}+\cvv{Q}-\cvv{R}$, then $\tilde{\cg}$
is the function over $\Sol$ that takes the value $1\in\mathbb{F}_2$ at all the three
points $P$, $Q$ and $R$, and vanishes otherwise. However,
$\mathbb{F}_2$ is now a field and $V(\Sol,\mathbb{F}_2)$ is simply a
vector space! Deciding whether $\tilde{\1}_I-\tilde{\1}_J$ belongs to
it is a simple matter 
requiring only linear algebra. 

Let us investigate this problem
further.
One way to characterize $V(\Sol,\mathbb{F}_2)$
as a subspace of $\F(\Sol,\mathbb{F}_2)$ is to compute equations of it. By
using the canonical scalar product, this reduces to computing 
$V(\Sol,\mathbb{F}_2)^\perp$ which means
the elements $\chi\in \F(\Sol,\mathbb{F}_2)$ such that
\begin{equation}
  \label{eq:4}
  \forall \cg=\cvv{P}+\cvv{Q}-\cvv{R}\in\D(\Sol),\quad
  \chi(P)+\chi(Q)=\chi(R)
\end{equation}
since any such $\chi$ verifies
\begin{equation}
  \label{eq:5}
  \forall g\in V(\Sol,\mathbb{F}_2),\quad\sum_{A\in\Sol}\chi(A)g(A)=0.
\end{equation}
We need a name for such elements of $V(\Sol,\mathbb{F}_2)^\perp$, and
we propose the name \emph{witness}. Let us start to do so on the english board. Let us determine a
function $\chi_0$. We first fix four values
on a square, for instance
\begin{figure}[!h]
  \centering
  \begin{minipage}[c]{0.4\textwidth}
     \centering \scalebox{0.3}{\includegraphics{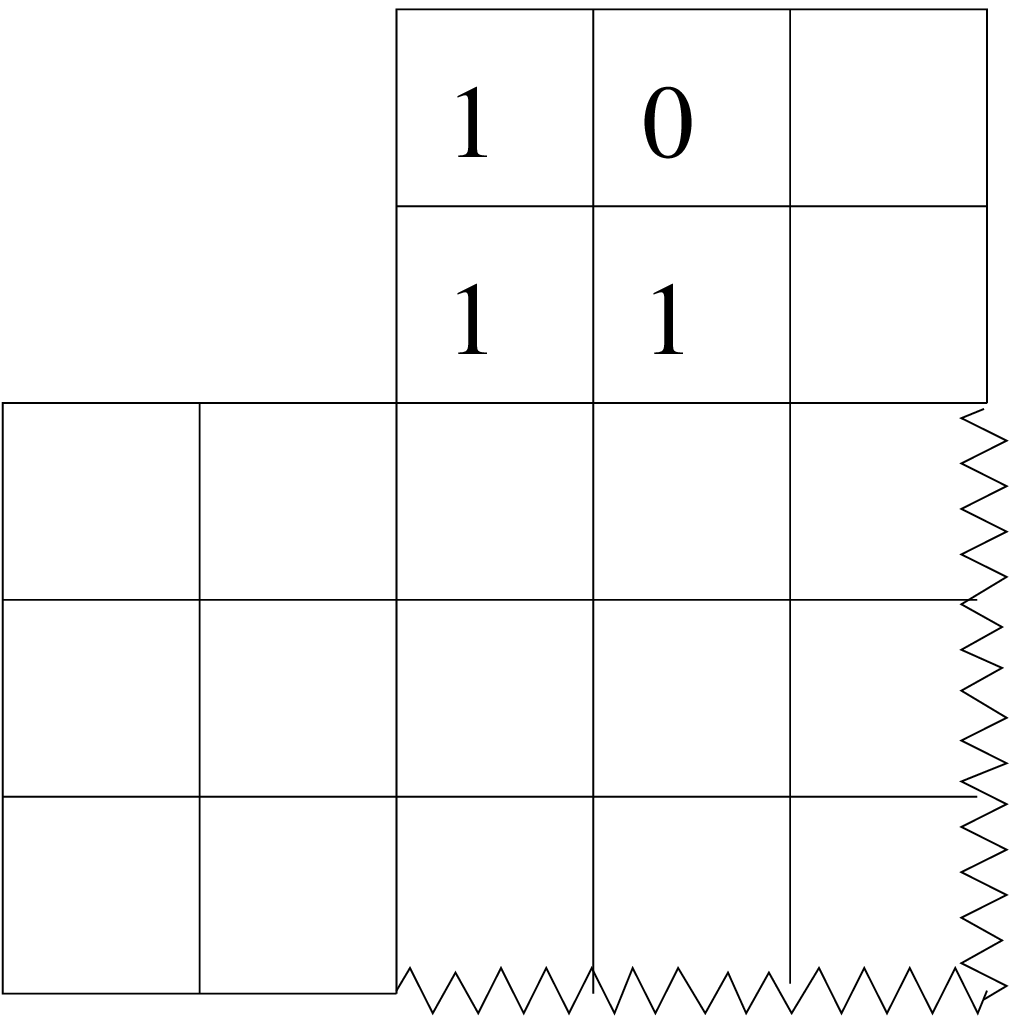}}%
     \caption{Starting values}
  \end{minipage}
\end{figure}
\FloatBarrier
By using \eqref{eq:4}, we can readily extend these values:
\begin{figure}[!h]
  \centering
  \begin{minipage}[c]{0.3\textwidth}
     \centering \scalebox{0.3}{\includegraphics{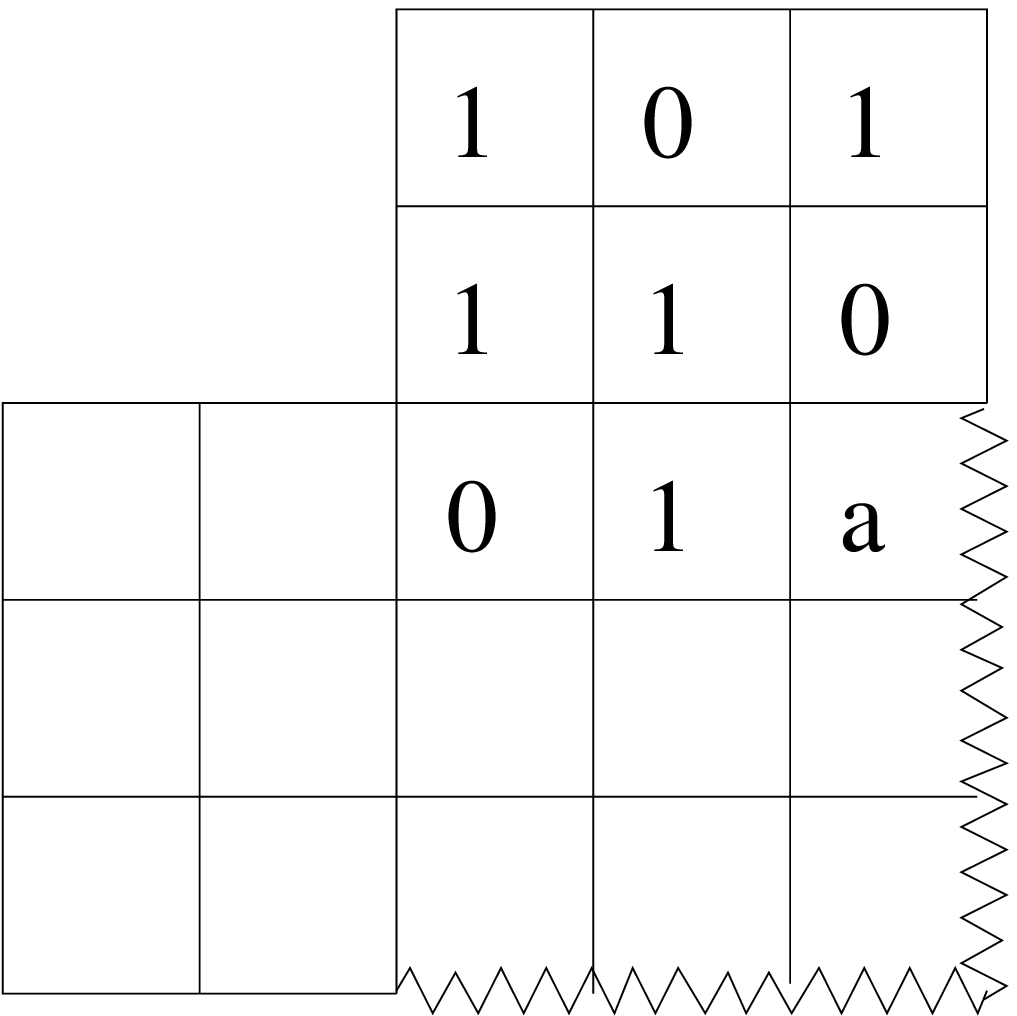}}%
     \caption{Extension}
  \end{minipage}
\end{figure}
\FloatBarrier
As it turns out, there are two ways to compute $a$: either by adding
the two values on the column above its square or the two on the line
containing it. The result is here the same $a=1$. We can use this
process to compute the values of $\chi_0$ on the full board. 

What is the dimension of $V(\Sol,\mathbb{F}_2)$ in this case? The
values on the initial square determine the values everywhere as we
have just now remarked, and there is thus 16 witnesses. But these
values are not linearly independant and  there are linearly generated by
the four
\begin{figure}[!h]
  \centering
  \begin{minipage}[c]{0.2\textwidth}
     \centering \scalebox{0.3}{\includegraphics{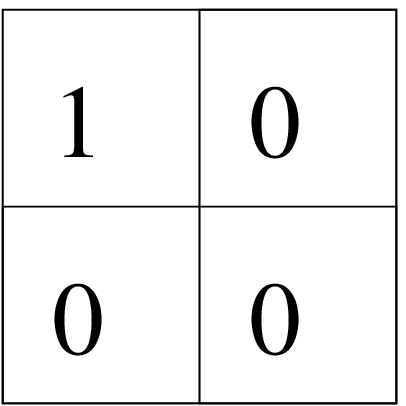}}%
  \end{minipage}
  \begin{minipage}[c]{0.2\textwidth}
     \centering \scalebox{0.3}{\includegraphics{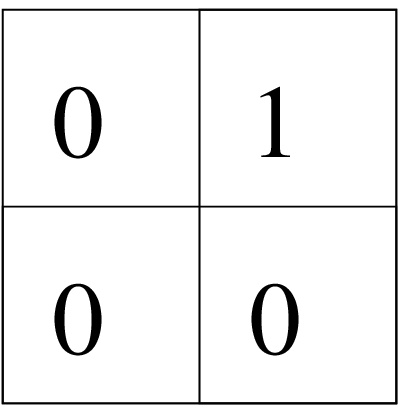}}%
  \end{minipage}
  \begin{minipage}[c]{0.2\textwidth}
     \centering \scalebox{0.3}{\includegraphics{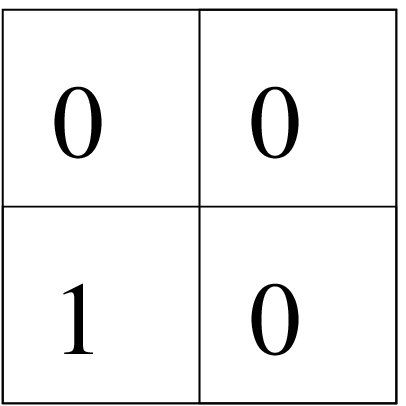}}%
  \end{minipage}
  \begin{minipage}[c]{0.2\textwidth}
     \centering \scalebox{0.3}{\includegraphics{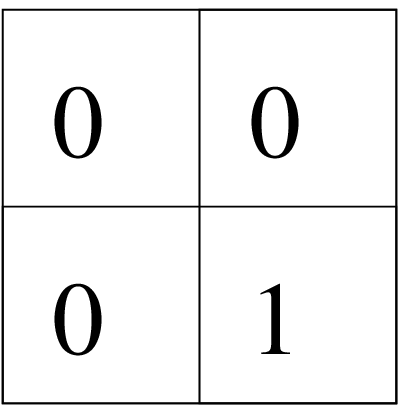}}%
  \end{minipage}
\end{figure}
\FloatBarrier
We can even use this process to extend the values to $\mathbb{Z}^2$.
This yields
\begin{figure}[!h]
  \centering
  \begin{minipage}[c]{0.8\textwidth}
     \centering \scalebox{0.3}{\includegraphics{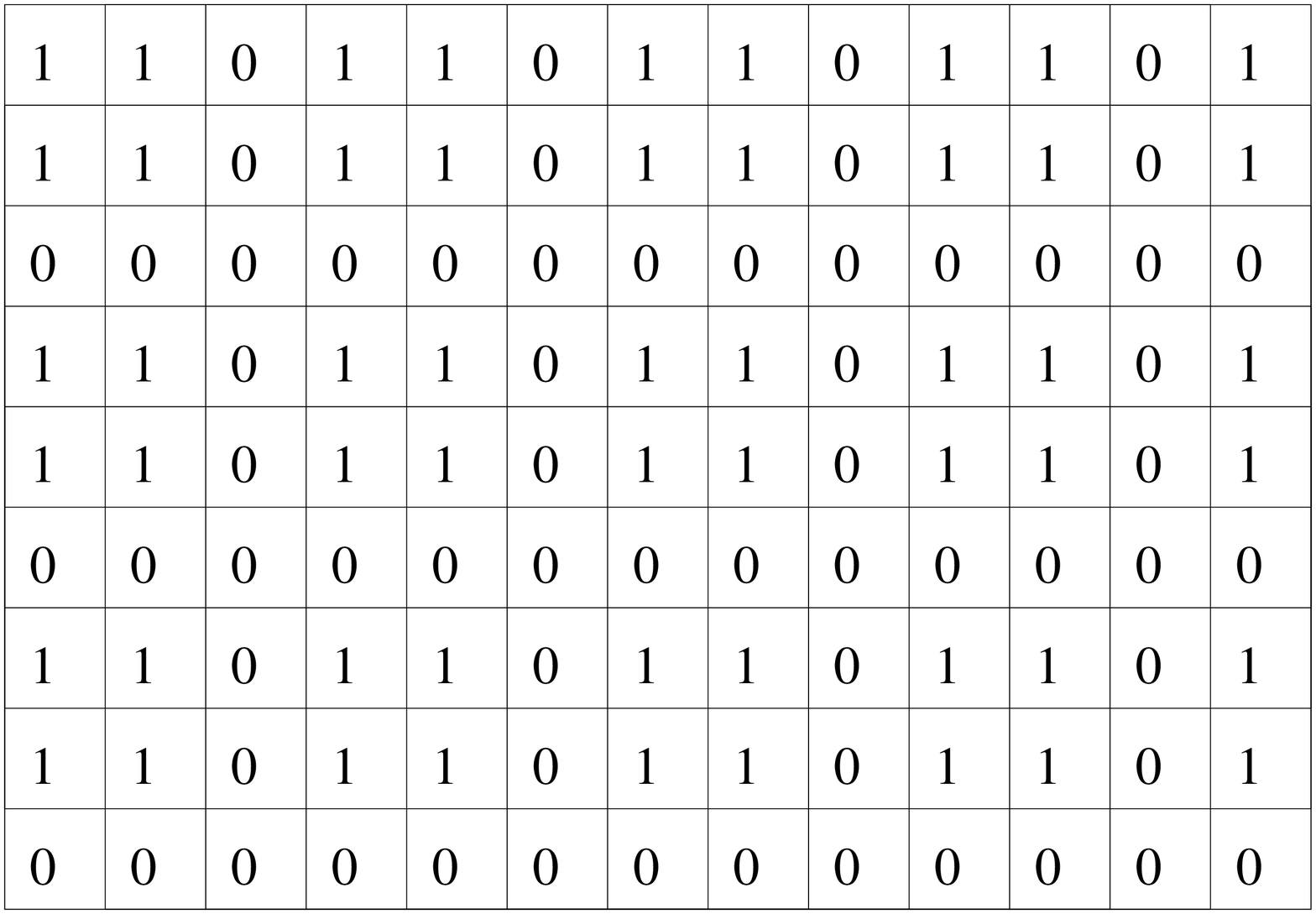}}%
     \caption{Over $\mathbb{Z}^2$}
  \end{minipage}
\end{figure}
\FloatBarrier
Now that the reader ses the regularity of this tiling, s.he will be
convinced that they can be extended to $\mathbb{Z}^2$. The way one
drops the english board on it yields for instance this witness:
\begin{figure}[!h]
  \centering
  \begin{minipage}[c]{0.8\textwidth}
     \centering \scalebox{0.3}{\includegraphics{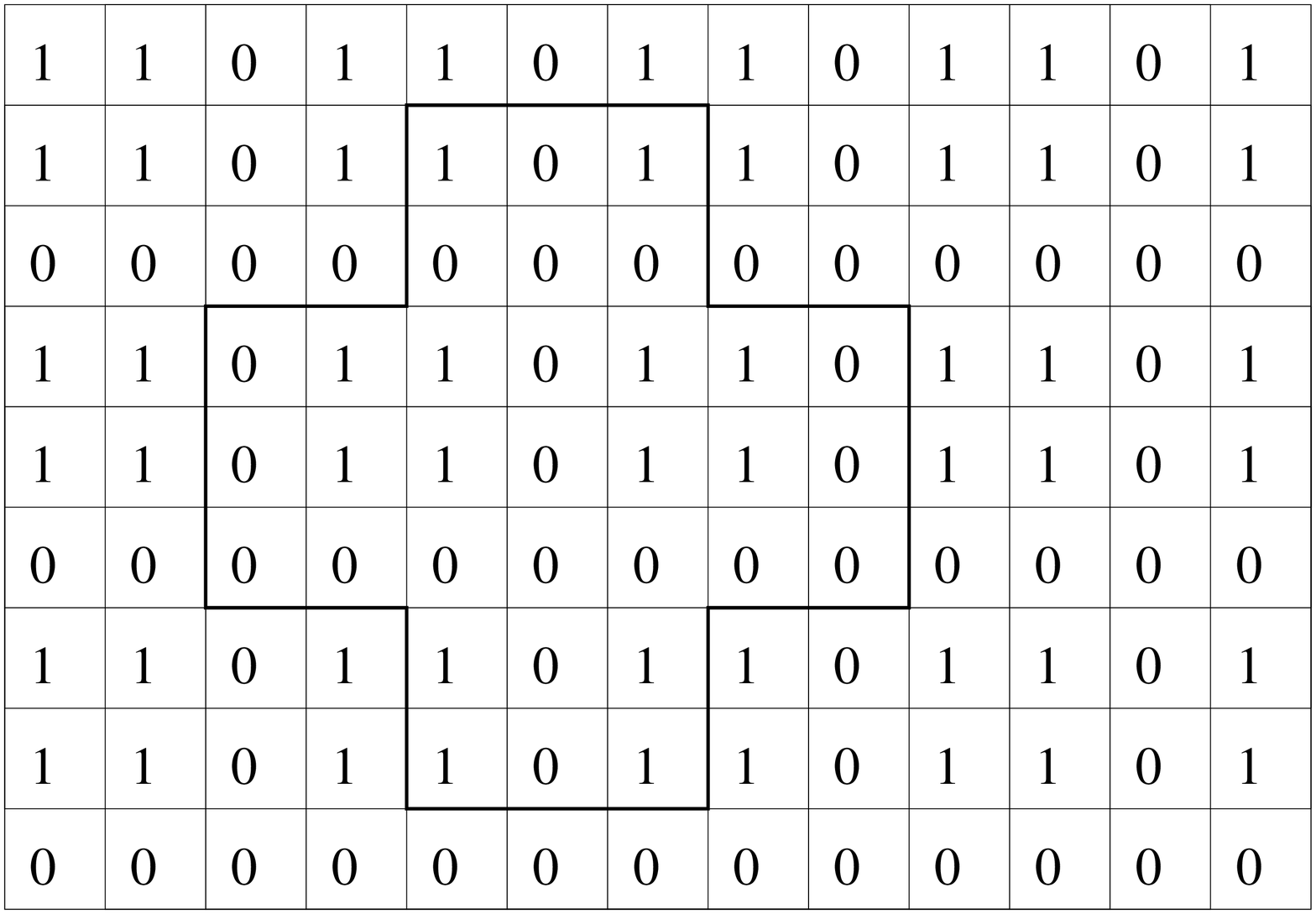}}%
  \end{minipage}
\end{figure}
\FloatBarrier
And once, the witnesses are determined, equations defining
$V(\Sol,\mathbb{F}_2)$ are obtained by taking the scalar product with
(a basis) of them.
A classical problem is to determine whether it is possible to start
with the french board filled with pegs, except for the central square
that is left empty and to end with only one peg. This can be shown to
be impossible by using the theory above, but we leave this pleasure to
the reader. 

This theory of witnesses is essentially what is called Reiss's theory~\cite{Reiss*57}, though it
is expressed with other words, and is present in Lucas's book~\cite{Lucas*91}. We say
"essentially" because they do not use any linear algebra and that
their way to reach this result is by using direct move together with
reversed ones (to undo a move). They obtain what they call
characteristic positions, which is equivalent to the equations
defining $V(\Sol,\mathbb{F}_2)$. This is however what is
presented~\cite{deBruijn*72}. There are still a distinction to be
made:
\begin{enumerate}
\item One can start from witnesses of $\mathbb{Z}^2$, restrict them to
  $\Sol$ and get witnesses for this board. This is called the
  rule-of-three. Of course, we get only a four independant equations
  that may not define  $V(\Sol,\mathbb{F}_2)$ fully. If the board is
  thick enough, for instance when there exists a defining square from
  which all the other values of the witnesses can be deduced, this is enough.
\item One can start from $V(\Sol,\mathbb{F}_2)$ and directly compute a
  basis of witnesses. This is required when the board is weakly
  connected (or even not connected!) and $V(\Sol,\mathbb{F}_2)^\perp$
  has dimension larger than~4. Several examples like that are given
  in~\cite{Avis-Deza-Onn*00}. 
\end{enumerate}

\section{The integer linear test and the lattice criterion}

Thinking bach in terms of $V(\Sol,\mathbb{Z})$, the lattice criterion
of  \cite{Avis-Deza-Onn*00} is to say that $\1_I-\1_J$ should belong
to $V(\Sol,\mathbb{Z})$. How is this test connected with the previous one?
Or, alternatively: we decided to reduce the problem modulo~2; Why not
try to do so modulo~3? Let us first note that we may identify 
$\F(\Sol,\mathbb{F}_2)$ with $\F(\Sol,\mathbb{Z})/2\cdot
\F(\Sol,\mathbb{Z})$ via
\begin{equation}
  \label{eq:19}
  \begin{array}{rcl}
    \tilde\ :\ \F(\Sol,\mathbb{Z}) &\rightarrow&\F(\Sol,\mathbb{Z})\\
    g&\mapsto&\tilde{g}\ :\
    \begin{array}[t]{rcl}
      \Sol&\rightarrow&\mathbb{F}_2\\
      P&\mapsto&g(P)\mod 2
    \end{array}
  \end{array}
\end{equation}
During this process, $V(\Sol,\mathbb{Z})$ is of course sent on
$V(\Sol,\mathbb{F}_2)$.
Let us state formally two questions we want to answer: 
\begin{enumerate}
\item Is $V(\Sol,\mathbb{Z})$ a lattice of full rank in
  $\F(\Sol,\mathbb{Z})$?
\item How to compute $\F(\Sol,\mathbb{Z})/V(\Sol,\mathbb{Z})$?
\end{enumerate}
In the sequel, we introduce a hypothesis on the geometry of the board
$\Sol$ that will enables us to answer fully these questions. It will
turn out that this will also exhibit the very tight link between the integer
linear test and the theory of witnesses, as exposed in the previous section.

If the two points $P$ and $R$ of $\Sol$ are extremities of a member of
$\D(\Sol)$, we say that $P$ and $R$ are \emph{neighbors} and we note
$P\voisin R$. The reflexive and transitive closure of this relation is an
equivalence relation, and if two points $A$ and $B$ are equivalent according
to it, we note $A\equiv B$. We can now state an important definition:
\begin{defi}
  A board $\Sol$ is said to be \emph{with no isolated point} if for every point $P$ of
  $\Sol$, there exists a point $Q\equiv P$ and which is the middle point of a move.
\end{defi}
Most boards will verify this hypothesis. It means that each $\equiv$-equivalence class
contains a middle point. However the number of such classes may vary. For a
sufficiently thick board, there will be exactly 4~classes, but there may be
more, if the board is not connected for instance, or contains thick chambers
very weakly connected by only one square. The reader will easily construct
examples of boards with no isolated point but were the number of classes is
larger than~4. 
The following Theorem is central in our discussion:
\begin{thm}
\label{th1}
If $\Sol$ is with no isolated point, then $2\F(\Sol,\mathbb{Z})\subset V(\Sol,\mathbb{Z})$.
\end{thm}

A final notation before sketching the proof:
if $\cg=\cvv{P}+\cvv{Q}-\cvv{R}\in\D(\Sol)$, we note $\cg'=-\cvv{P}+\cvv{Q}+\cvv{R}$
the reversed move (with equal middle point).
\begin{dem}
  We show that for every $P\in\Sol$, we have $2\cvv{P}\in V(\Sol,\mathbb{Z})$.
  If $P$ is a middle point, say of the move $\cg$, then $2\cvv{P}=\cg+\cg'$ belongs
  to $V(\Sol,\mathbb{Z})$.  Otherwise, there exists a chain $P=P_0\voisin P_1\voisin \dots\voisin
  P_n$ where $P_n$ is a middle point. Furthermore, by definition, there exists
  $\cg_i\in\D(\Sol)$ such that $2\cvv{P}_i-2\cvv{P}_{i+1}=\cg_i-\cg'_i$ for every
   $i=0,\dots,n-1$. Finally, we can also write $2\cvv{P}_n=\cg_n+\cg'_n$  for some
  $\cg_n\in\D(\Sol)$. Summing up all these equations, we reach
  \begin{equation*}
    2\cvv{P}_0=
    \cg_0-\cg'_0+\cg_1-\cg'_1+\dots+\cg_{n-1}-\cg'_{n-1}+\cg_n+\cg'_n\in V(\Sol,\mathbb{Z}),
  \end{equation*}
  which is the required conclusion since  $P=P_0$.
\end{dem}
This Theorem has several consequences. First of all, on such boards, the
$\mathbb{Q}$-vector spanned by the $\cg$'s (that would be
$V(\Sol,\mathbb{Q})$) is the whole space: $V(\Sol,\mathbb{Z})$ is a sublattice
of $\F(\Sol,\mathbb{Z})$ of full rank. Let us note the following Lemma
that will be required later:
\begin{lem}\label{count}
  If $\Sol$ has no isolated points, we have $|\Sol|\le |\D(\Sol)|\le 4|\Sol|-8$.
\end{lem}
\begin{dem}
  The lower bound comes from the fact that $\D(\Sol)$ generates
  $\F(\Sol,\mathbb{Q})$. For the upper bound, count horizontal and vertical
  moves separately. For the horizontal (resp. vertical) ones, count the moves according to their
  left-hand side (resp. lower) point. The lemma follows readily. 
\end{dem}

As a main consequence, we have the following Theorem.
\begin{thm}\label{struct}
  Assume $\Sol$ to be with no isolated point and let
  $g\in\F(\Sol,\mathbb{Z})$. Then 
  \begin{equation*}
    g\in V(\Sol,\mathbb{Z})
    \iff
    \tilde{g}\in V(\Sol,\mathbb{F}_2).
  \end{equation*}
  (See \eqref{eq:19} for the definition of $\tilde{g}$).
\end{thm}
\begin{dem}
  Indeed, the direct implication is obvious, while the reversed one
  follows from Theorem~\ref{th1}: we know that
  $g\in V(\Sol,\mathbb{Z})+2\cdot\F(\Sol,\mathbb{Z})$
  but this last space is nothing but $V(\Sol,\mathbb{Z})$.
\end{dem}

This Theorem tells us that the lattice criterion is \emph{not}
stronger than Reiss's theory, when properly understood, and provided
we restrict our attention to non-pathological boards.
In fact \cite{Avis-Deza-Onn*00} do not even give a single example when 
reduction modulo~2 does not solve the problem. Here is one:
\begin{figure}[!h]
  \centering
  \begin{minipage}[c]{0.4\textwidth}
     \centering \scalebox{0.3}{\includegraphics{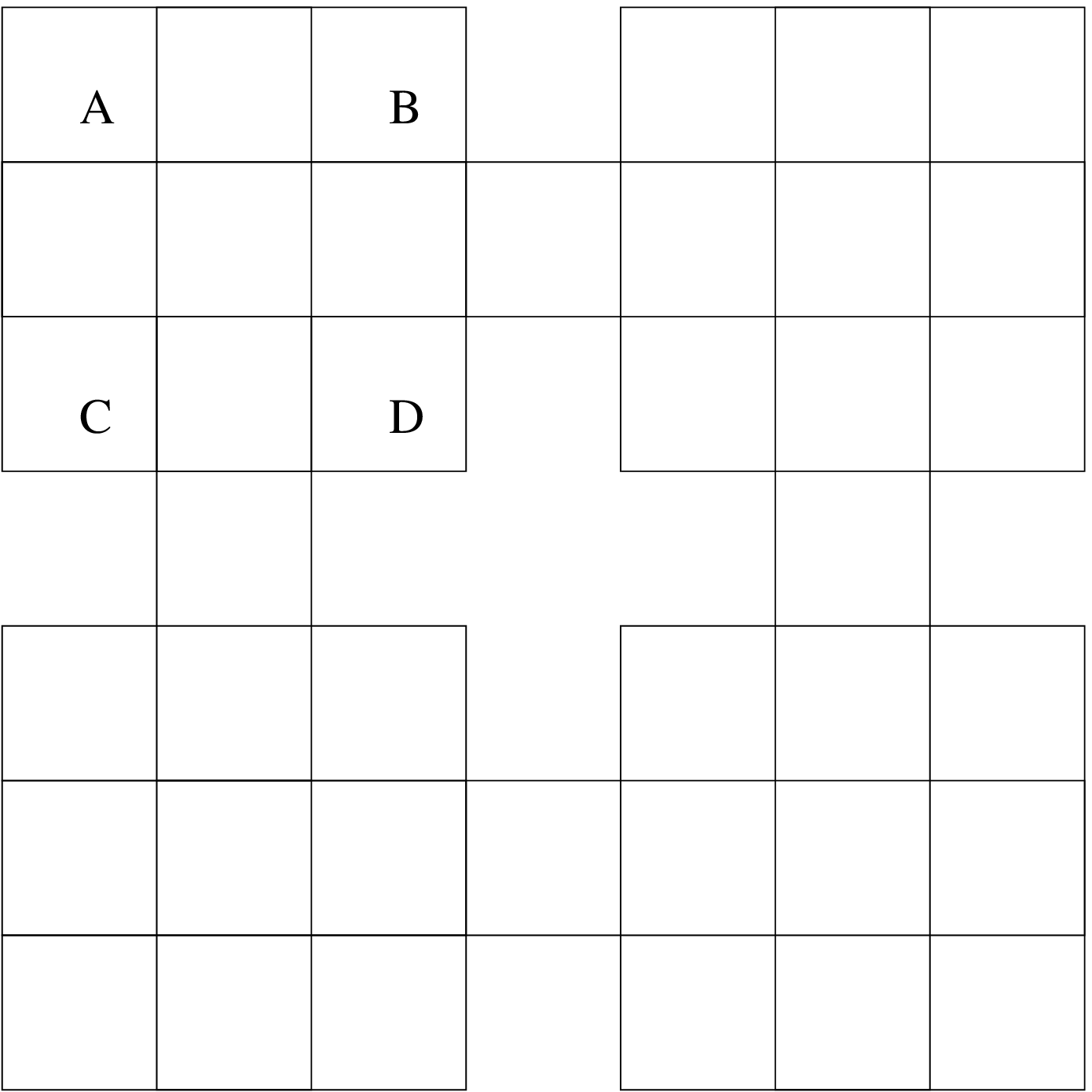}}%
     \caption{A pathological board}\label{patho}
  \end{minipage}
\end{figure}
\FloatBarrier
The total number of pegs on the squares $A$, $B$, $C$ and $D$ remains
constant. It is not difficult to see that this example is in fact
general and we have:
\begin{thm}\label{restruct}
  A board $\Sol$ is with no isolated point if and only if
  $V(\Sol,\mathbb{Z})$ has maximal rank in $\F(\Sol,\mathbb{Z})$. 
\end{thm}
On boards with no isolated points, reducing the situation modulo any
odd integer is not going to give any information; indeed
Theorem~\ref{struct} implies (after some work) that 
\begin{equation*}
  V(\Sol,\mathbb{Z}/m\mathbb{Z})=\F(\Sol,\mathbb{Z}/m
  \mathbb{Z})
  \quad(\text{whenever $m$ is odd}).
\end{equation*}
Notice finally that $\F(\Sol,\mathbb{Z})/V(\Sol,\mathbb{Z})$
is simply a product of copies of $\mathbb{Z}/2\mathbb{Z}$ in this case.
It is not difficult to tackle the case with isolated points by
generalising the reasoning used for the board drawn figure~\ref{patho}, and get that 
 $\F(\Sol,\mathbb{Z})/V(\Sol,\mathbb{Z})$ is always a product of
 copies of $\mathbb{Z}/2\mathbb{Z}$ with copies of $\mathbb{Z}$. These
 results have no influence on what we develop hereafter, so we do not
 provide any formal proof.

\section{Resource counts, pagoda functions and the linear test in
  non-negative rationals}

The next main step takes place in 1961/1962 at Cambridge university
when J.H.~Conway led a group of students (among which were Beasley)
that studied this game. They came out with another and different test,
also clearly explained in~\cite{Berlekamp-Conway-Guy*82} and that we
now describe.

This test exploits the fact that~\eqref{eq:1} has non-negative
coefficients, i.e. the test consists in writing that, if we can go from $I$ to
$J$ with legal moves, then
\begin{equation}
  \label{eq:6}
  \1_I-\1_J\in V^+(\Sol,\mathbb{Q}).
\end{equation}
As it turns out, $V^+(\Sol,\mathbb{Q})$ is a cone in a vector space, and
determining whether a point belongs to it or not is fast. We know generators
of this cone (the elements of $\D(\Sol)$; they can be shown to be generator of
its extreme half-lines), and it would be interesting to determine equations for its
facets. The paper~\cite{Avis-Deza*01} gives properties of these facets.
In~\cite{Berlekamp-Conway-Guy*82} as well as in~\cite{Beasley*92}, so called
\emph{resource counts} or \emph{pagoda functions} are introduced. These are
functions $\pi$ on $\Sol$ such that
\begin{equation}
  \label{eq:7}
  \forall \cg=\cvv{P}+\cvv{Q}-\cvv{R}\in\D(\Sol),\quad
  \pi(P)+\pi(Q)\ge\pi(R).
\end{equation}
As a consequence, for any such function and if $g$ belongs to
$V^+(\Sol,\mathbb{Q})$, one has 
\begin{equation}
  \label{eq:8}
  \langle \pi,g\rangle =\sum_{A\in\Sol}\pi(A)g(A)\ge 0.
\end{equation}
In particular, if one can derive $J$ from $I$ with legal moves, then $ \langle
\pi,\1_I\rangle$ is not less than $ \langle \pi,\1_J\rangle$. Here are some examples

\begin{figure}[!h]
  \centering
  \begin{minipage}[c]{0.26\textwidth}
     \centering \scalebox{0.20}{\includegraphics{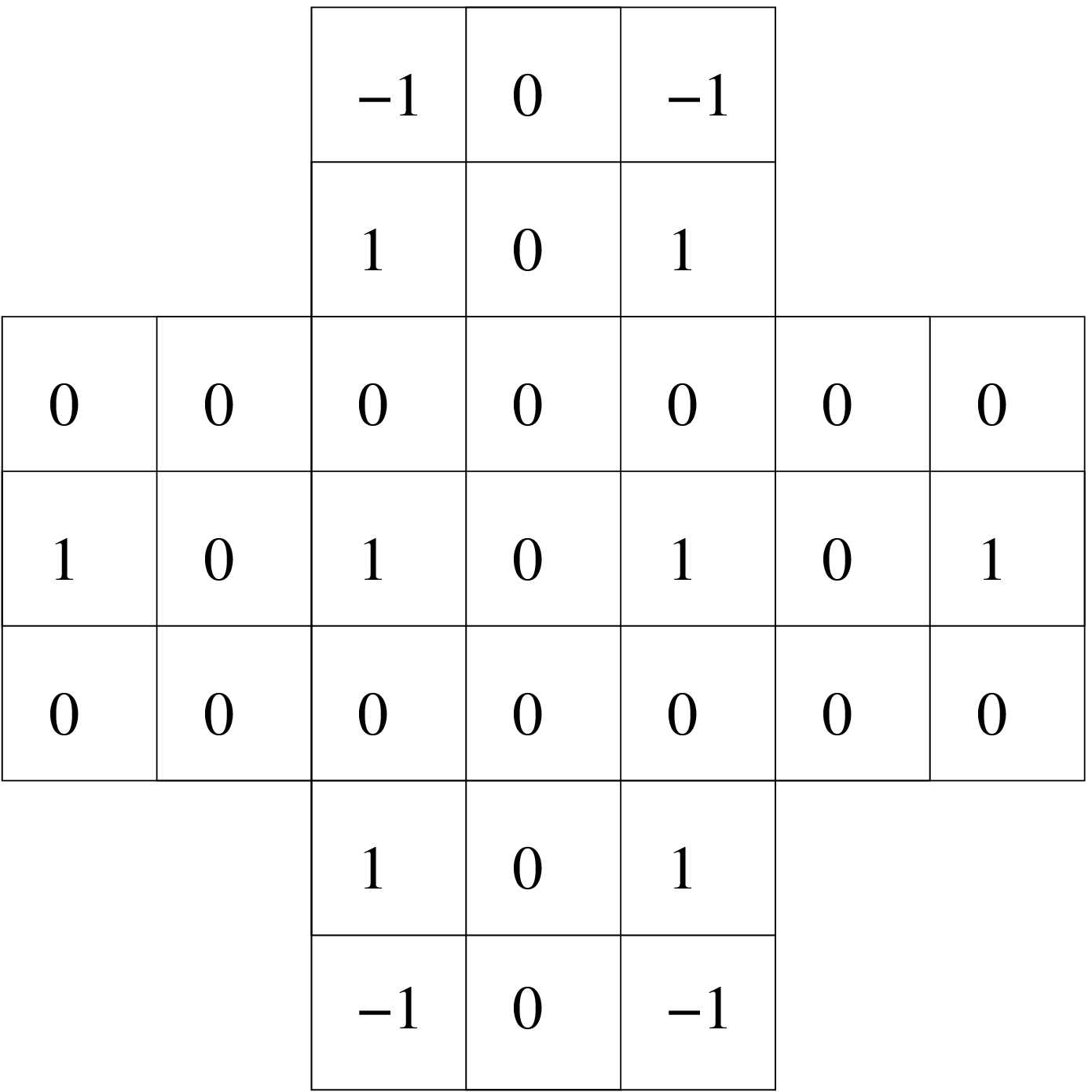}}%
     \caption{A resource count}
  \end{minipage}\quad\quad
  \begin{minipage}[c]{0.26\textwidth}
     \centering \scalebox{0.20}{\includegraphics{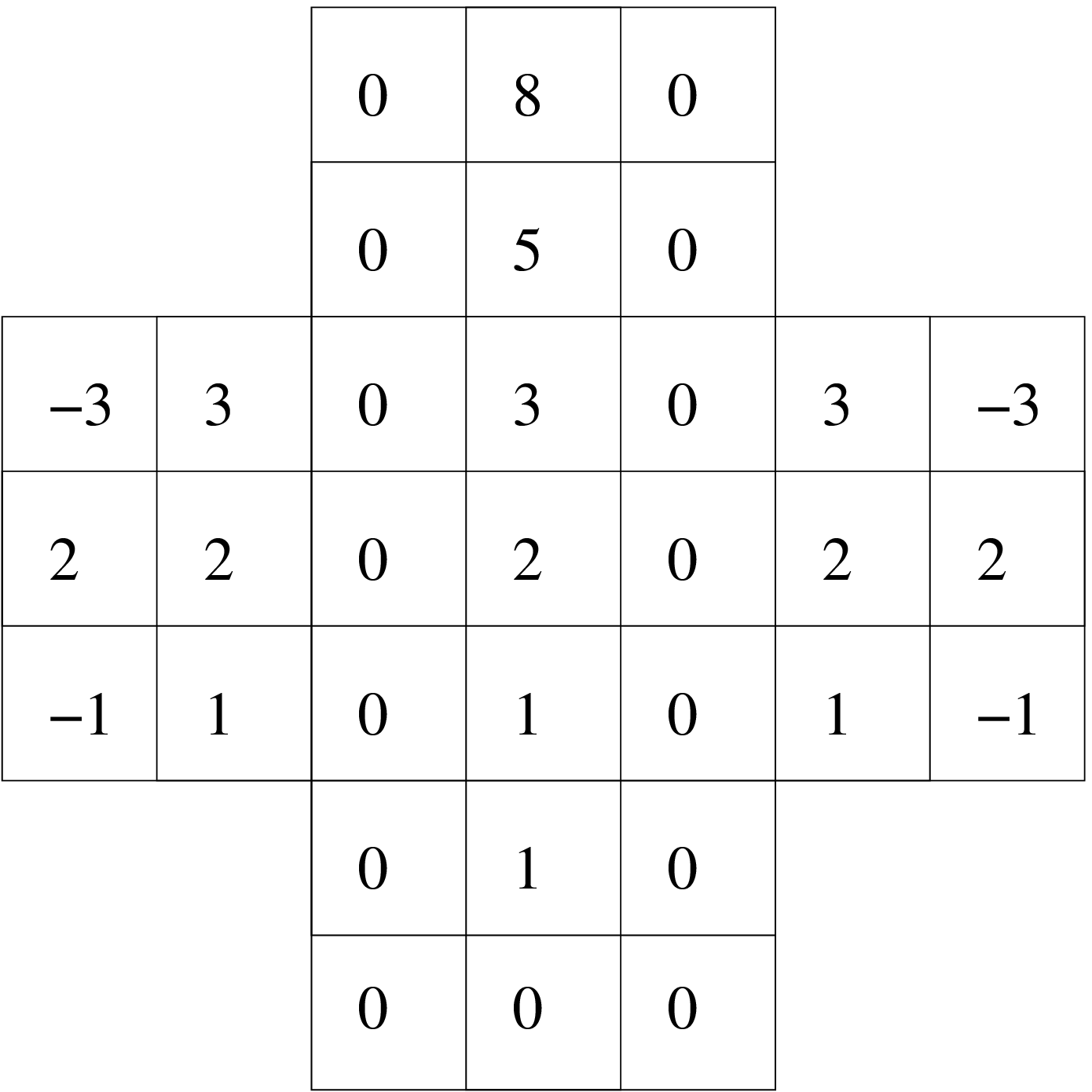}}%
     \caption{Another resource count}
  \end{minipage}\quad\quad
  \begin{minipage}[c]{0.26\textwidth}
     \centering \scalebox{0.20}{\includegraphics{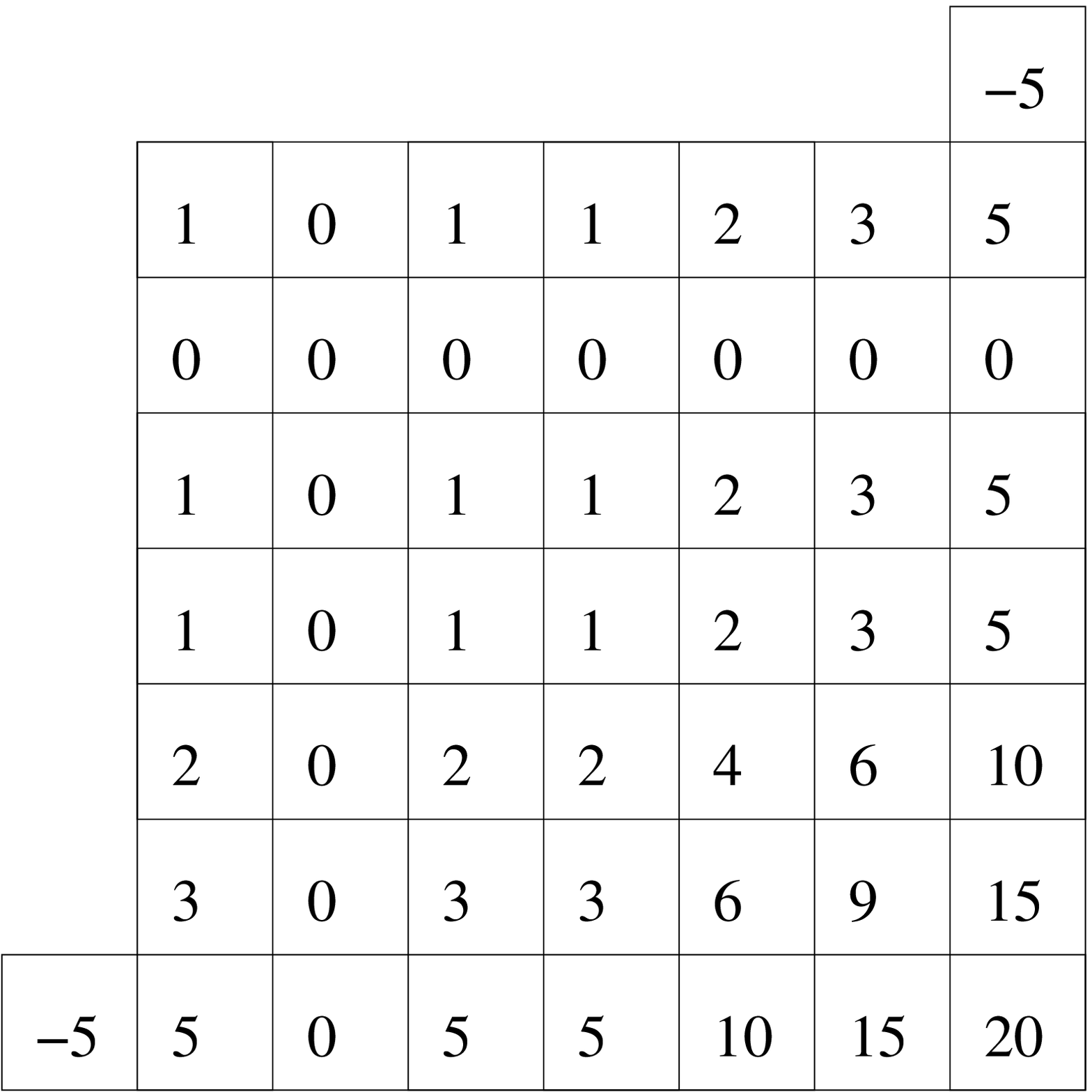}}%
     \caption{A third resource count}
  \end{minipage}
\end{figure}

\FloatBarrier

Determining which of these corresponds to equations of facets would be very
valuable, but their structure seems too intricate to classify
them in a small number of regular families. For instance , a direct
computation in case of the english board stumbles on the fact
that there are an enormous quantity of such facets for a human eye to be able
to look at them and derive some patterns. It is not sure that
this path is blocked, though I tend to believe it is. 

We do not dwell any
further in this part of the theory since it is extremely well exposed and
detailed in~\cite{Berlekamp-Conway-Guy*82},~\cite{Beasley*92} and on a number
of web pages. The reader will most probably better unterstand the
strength of this theory by looking at section~\ref{thickness} of this paper.

We should stress out here that the approach of this Cambridge group is
commonly reduced to the use of real-valued "pagoda" functions as
above. This is
an extremely minimal understanding of their work and for instance does
not account for the GNP balance sheet, what Beasley in~\cite{Beasley*92} calls
Conway's balance sheet in his chapter~6; this one is however one of
the main tool of~\cite{Berlekamp-Conway-Guy*82}.  It mixes \emph{integer
valued} pagoda functions together with such functions with values in
$\mathbb{F}_2$.  Beasley's use of pagoda functions which he calls
ressource counts (see chapter~5 of~\cite{Beasley*92}) relies already
on the integer character of the values taken: that is how he builds
his "move map".

The GNP diagram, or GNP balance sheet, is somewhat off our
framework, and is in fact superseded by the next test.

\section{The linear test in non-negative integers}

\begin{figure}[!h]
  \centering
  \begin{minipage}[c]{0.3\textwidth}
     \centering \scalebox{0.3}{\includegraphics{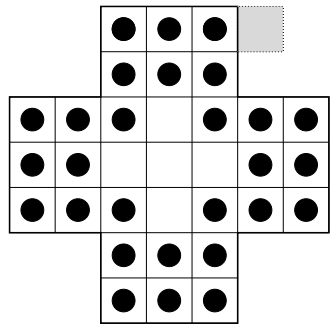}}%
     \caption{Impossible}\label{ce11}
  \end{minipage}
\end{figure}

The third test consists in combining both preceding ideas and write that
 if we can go from $I$ to
$J$ with legal moves, then
\begin{equation}
  \label{eq:9}
  \1_I-\1_J\in V^+(\Sol,\mathbb{Z}).
\end{equation}
This time, deciding that an element belongs to the integer points of a cone is
NP-hard, but in practice, it takes only some fraction of a second on an
english board (this was not the case in 1962!).
We have of course
\begin{equation}
  \label{inter}
  V^+(\Sol,\mathbb{Z})\subset V(\Sol,\mathbb{Z})\bigcap V^+(\Sol,\mathbb{Q})
\end{equation}
and this inclusion is strict, even when one restricts our attention to
differences of characteristic functions. For instance this test shows that one
cannot go from the position of figure~\ref{ce11} to only the central peg
while the rational and integer linear tests are passed. 
This example is interesting in showing the impact of the board, for it is
feasible in legal moves if we add to the english board the grey square on the
upper right side.

\begin{figure}[!h]
  \centering
  \begin{minipage}[c]{0.4\textwidth}
     \centering \scalebox{0.4}{\includegraphics{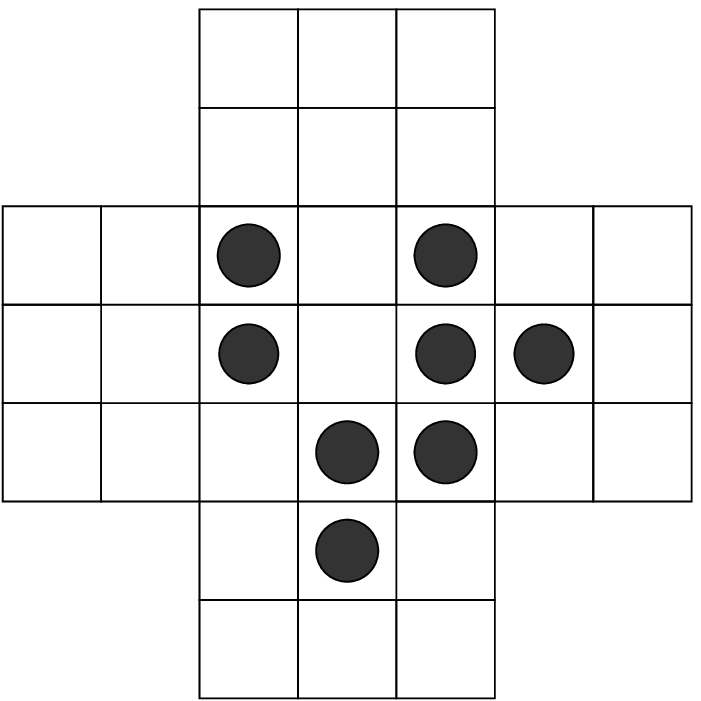}}%
     \caption{Starting position}\label{ce9}
  \end{minipage}\quad\quad
  \begin{minipage}[c]{0.4\textwidth}
     \centering \scalebox{0.4}{\includegraphics{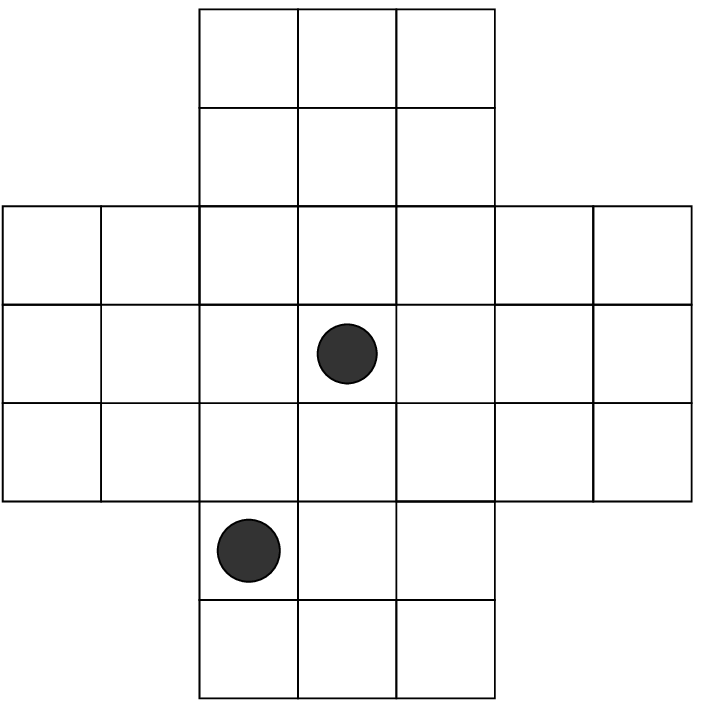}}%
     \caption{Ending position}\label{ce10}
  \end{minipage}
\end{figure}

We present a smaller counterexample in figures~\ref{ce9}
and~\ref{ce10} that enables easier direct computations.

When $g$ belongs to this intersection (i.e. the right-hand side
of~\eqref{inter}) the denominators in a non-negative 
writing do not seem to be any worse than $1/2$. Here is the conjecture we make:
\begin{conj}
  If $\Sol$ has no isolated points, then
  \begin{equation*}
    V(\Sol,\mathbb{Z})\bigcap V^+(\Sol,\mathbb{Q})\subset\tfrac12V^+(\Sol,\mathbb{Z}).
  \end{equation*}
\end{conj}
Here is another related conjecture that may be easier to handle (and maybe
easier to disprove!).
\begin{conj}
Let $\mathcal{B}\subset \D(\Sol)$ be a basis of $V(\Sol,\mathbb{Q})$.
If $\Sol$ has no isolated points, then
\begin{equation*}
2\F(\Sol,\mathbb{Z})\subset
V(\mathcal{B})=\sum_{\cg\in\mathcal{B}}\mathbb{Z}\cdot \cg.
\end{equation*}
\end{conj}
The condition on $\Sol$ cannot be removed since  it
is equivalent to $V(\Sol,\mathbb{Z})$ being of full rank.

At this point, we have described the situation and we hope the reader is now
able to understand properly what is what. The theory so far has two drawbacks:
it draws only on properties of $\1_I-\1_J$, and it does not use the order in
which the moves are played. Our next criteria, the simple quadratic test, will
not go beyond this abelian nature, but will break the first hurdle.
It is better to investigate the game a bit further before exposing it.


\section{How integer linear programming is used}

The cone $V^+(\Sol,\mathbb{Z})$ is determined by the set $\D(\Sol)$ of
generators. Let us introduce the notation $\cvv{\cg}$ for the function over
$\D(\Sol)$ that is 1 in $\cg$ and 0 everywhere else. We consider the
map
\begin{equation}
  \begin{array}{rccl}
    \Psi :& \F(\D(\Sol),\mathbb{Z})&\rightarrow&V(\Sol,\mathbb{Z})\\[1em]
    &\displaystyle F=\sum_{\cg\in\D(\Sol)}x(\cg)\,\cvv{\cg}
    &\mapsto& \displaystyle\sum_{\cg\in\D(\Sol)}x(\cg)\,\cg.
  \end{array}
\end{equation}
The integer linear program we write is simply to minimize any linear form of
the $(x(\cg))_\cg$ subject to the constraints
\begin{equation*}
  \forall\cg\in\D(\Sol), \ x(\cg)\ge 0,\quad\text{and}\quad\Psi(F)=\1_I-\1_J.
\end{equation*}
The linear form we choose is usually $\sum_\cg x(\cg)$ since we know what
should be its value if a solution exists.

\section{Thickness of a move}
\label{thickness}
Given a problem, say from $I$ to $J$, we define the \emph{thickness} of the move
$\cg$ to be the maximum number of times this move can be used, whatever sequence of
legal moves $\cg_1,\cg_2,\dots, \cg_k$ we choose. This thickness is zero
allover if the problem is not feasible. 
In general, given $h\in V^+(\Sol,\mathbb{Z})$, we shall speak of the
\emph{thickness of $\cg$ at $h$}.
Computing this quantity is naturally
difficult, but we can bound it from above and even provide a uniform bound
for it. The main Theorem reads as follows
\begin{thm}
  Let $h\in V^+(\Sol,\mathbb{Z})$, $\cg_0\in\D(\Sol)$ and $\pi$ be a
  resource count on $\Sol$ such that $\langle\pi,{\cg_0}\rangle=1$.
  The move $\cg_0$ can appear at most $\langle{\pi},{h}\rangle$ in any
  writing of $h$ as a linear combination of elements of $\D(\Sol)$
  with non-negative integer coefficients.
\end{thm}
The scalar product $\langle{\pi},{h}\rangle$ is defined in
\eqref{eq:8}. We can derive absolute bounds from this Theorem by using
a variant of a resource count already used by Conway. First note that
we are interested only in the case $h=\1_I-\1-J$ which implies that
$|h(A)|\le 1$ for all $A\in\Sol$.  Now let $\rho=(\sqrt{5}-1)/2$ be a
solution of $x^2+x=1$. To each point $(a,b)\in\mathbb{Z}^2$, we
associate the weight $\pi(a,b)=\rho^{|a|+|b|}$. Next, we drop our
board $\Sol$ on $\mathbb{Z}^2$ in such a way that the middle point of
$\cg_0$ be the $(0,0)$ element. The reader will check that the
restriction of $\pi$ to $\Sol$ is a resource count on $\Sol$ which we
denote again by $\pi$. We have $\langle\pi,{\cg_0}\rangle =1$, while
\begin{equation*}
  |\langle{\pi},{h}\rangle|\le \langle{\pi},\1_{\Sol}\rangle\le 8\rho+13=17.944\cdots.
\end{equation*}
This short argument show that the thickness of any move on any board is
bounded above by $17$. This is most probably a way too large majorant
(reaching a thickness of $4$ is already extremely difficult, and it can be
shown on using better resource counts that the maximal thickness on the english
board is at most~$5$), but it is \emph{universal}, i.e. independant of
the board we choose.

A similar argument is also the main ingredient of~\cite{Ravikumar*97} (see
Theorem~3.1 therein, with most probably a wrong computation at the end. The~26
of this result is to be replaced by a~34 but this leaves the rest of the
argument intact), and is the basis on which rely the low complexity results. 

Given a problem, we can refine this upper bound by selecting a more
appropriate resource count. Furthermore, once a majorant is given, say $m$, we
can check whether $\1_I-\1_J-m\cg_0$ is feasible or not (this means, whether
it passes whichever test we select). If not, we decrement $m$ and repeat the process.

\section{A simple quadratic test}

Let us consider the two following problems: we are to go from the left
hand side position with only the black pegs (or with the grey peg added) to
the right hand side one with 
a sole black peg (or with the grey peg added). Both problems pass the positive integer
test. 
\begin{figure}[!h]
  \centering
  \begin{minipage}[c]{0.4\textwidth}
     \centering \scalebox{0.3}{\includegraphics{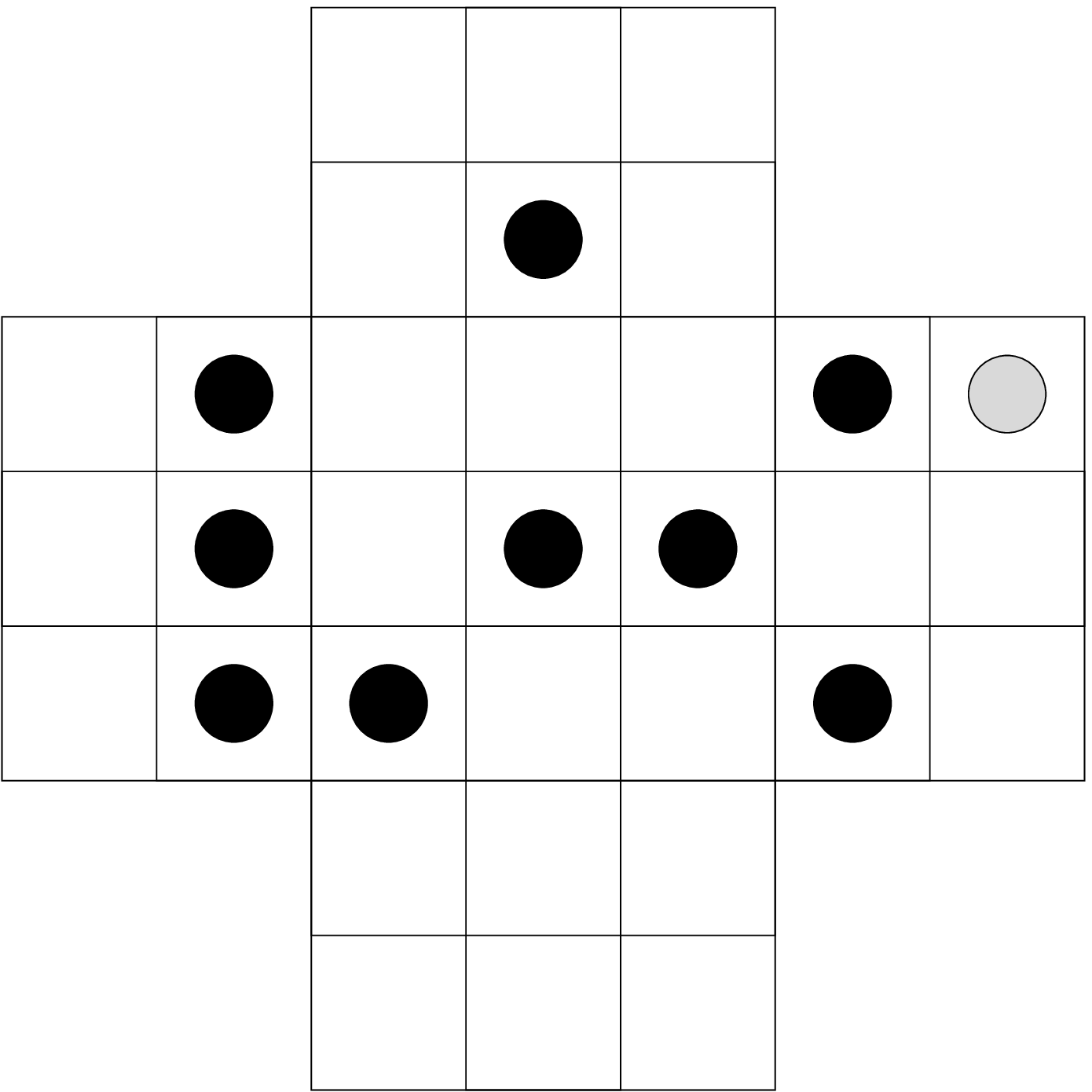}}%
     \caption{Starting position\label{cvvd1}}
  \end{minipage}\quad\quad
  \begin{minipage}[c]{0.4\textwidth}
     \centering \scalebox{0.3}{\includegraphics{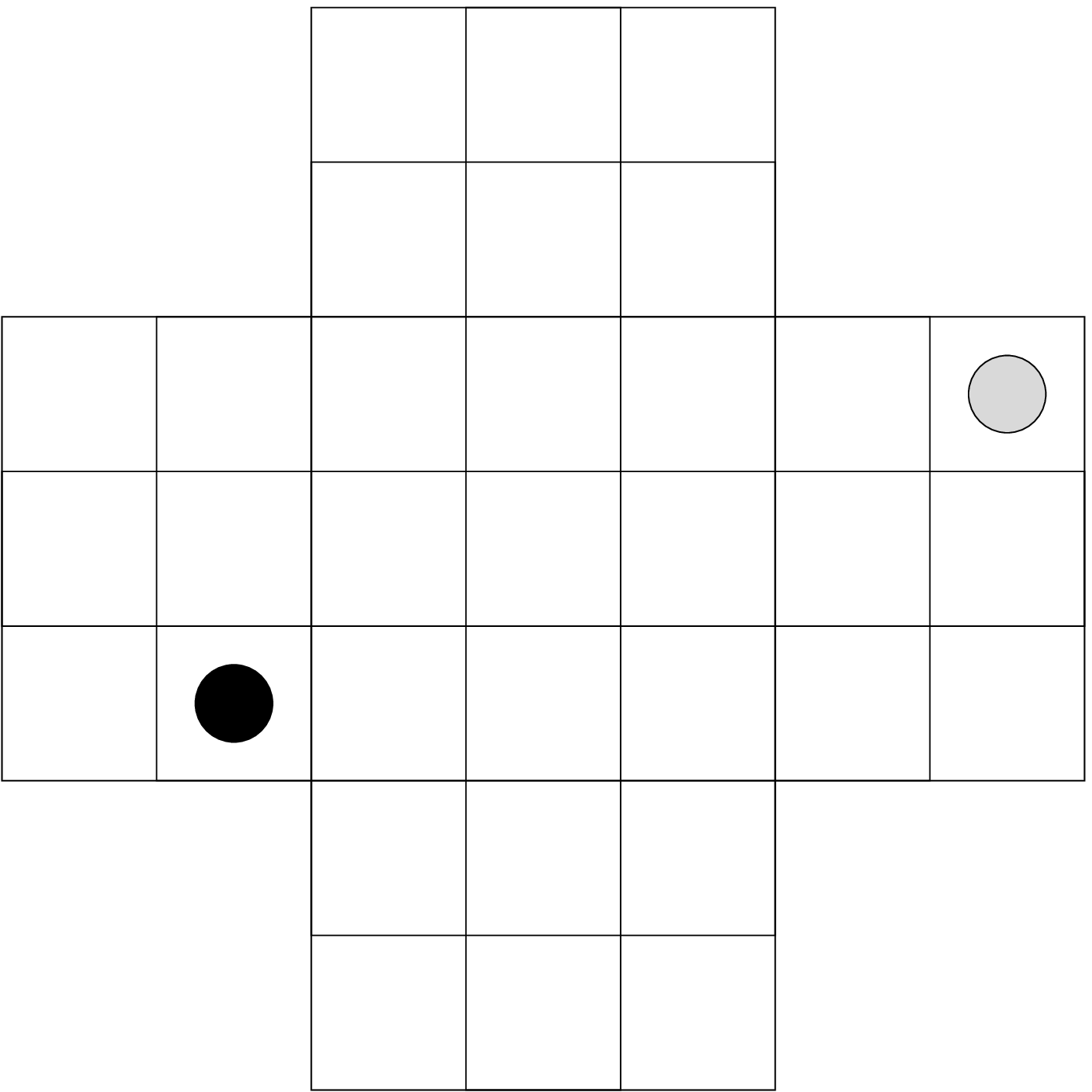}}%
     \caption{Ending position\label{cvvf1}}
  \end{minipage}
\end{figure}
The reader will easily check that the larger problem (with the grey
peg) is in fact doable in \emph{legal moves}, which implies that
\emph{no} test relying only on $\1_I-\1_J$ would be able to 
show the first problem to be impossible.
The quadratic test we propose now is however able to show this impossibility.
\FloatBarrier

Let us start our description of the quadratic test.

\subsection{The geometrical support}
To each couple $(A,B)\in\Sol\times\Sol$, we associate
a symbol $A\sqc B$, to which we add the property
\begin{equation}
  A\sqc B=B\sqc A.
\end{equation}
We set
\begin{equation}
  \Sol\sqc\Sol=\bigl\{A\sqc B,A,B\in\Sol\bigr\}.
\end{equation}
We next consider functions on $\Sol\sqc\Sol$. We denote by 
$\cvv{A\sqc B}$ the function that is 1 on $A\sqc B$ and 0 everywhere else.
Note that $\cvv{A\sqc B}=\cvv{B\sqc A}$. 
We go from $\F(\Sol,\mathbb{Q})^2$ to $\F(\Sol\sqc\Sol,\mathbb{Q})$ by 
\begin{eqnarray*}
  \sqc\quad:\quad
  \F(\Sol,\mathbb{Q})\times\F(\Sol,\mathbb{Q})&\rightarrow&\
  \F(\Sol\sqc\Sol,\mathbb{Q})\\ 
  (g_1,g_2)&\mapsto&\displaystyle
  g_1\sqc g_2 = \sum_{A,B\in\Sol}g_1(A)g_2(B)\cvv{A\sqc B}.
\end{eqnarray*}
Notice that the value of $g_1\sqc g_2$ on $\cvv{A\sqc B}$ is
$g_1(A)g_2(B)+g_1(B)g_2(A)$ if $A\neq B$ and $g_1(A)g_2(A)$ if $A=B$.

\subsection{The effect of legal moves}
Assume now that we can go from $I$ to $J$ by the legal move
$\cg\in\D(\Sol)$. We have
\begin{equation*}
  1_I\sqc\1_I=(\1_J+\cg)\sqc(\1_J+\cg)
  =\1_J\sqc\1_J + \cg\sqc\1_J+\1_J\sqc \cg+\cg\sqc \cg.
\end{equation*}
On using the identity $\cg\sqc\1_J=\1_J\sqc \cg$, we reach
\begin{equation*}
  \1_I\sqc\1_I=\1_J\sqc\1_J+(2\1_J+\cg)\sqc \cg.
\end{equation*}
We note that
\begin{equation*}
  2\1_J+\cg=\sum_{\substack{A\in J\\ \cg(A)=0}}2\cvv{A}
  +|\cg|,
\end{equation*}
from which we infer
\begin{equation}
  \label{eq:10}
  \1_I\sqc\1_I=\1_J\sqc\1_J+|\cg|\sqc \cg+\sum_{\substack{A\in J\\ \cg(A)=0}}2\cvv{A}\sqc\cg.
\end{equation}
This is the equation we want to exploit; we do so in pretty much the same
way we exploited~\eqref{eq:1}. We set
\begin{multline}
  \D(\Sol\sqc\Sol)=\{2\cvv{A}\sqc \cg,\quad  A\in\Sol, 
  \cg\in\D(\Sol)/\cg(A)=0\}
  \\\bigcup\
  \{|\cg|\sqc \cg,\quad  \cg\in\D(\Sol)\}.
\end{multline}
Note that if $\cg=\cvv{P}+\cvv{Q}-\cvv{R}$ then 
\begin{equation}
  |\cg|\sqc \cg=\cvv{P\sqc P}+2\cvv{P\sqc Q}+\cvv{Q\sqc Q}-\cvv{R\sqc R}.
\end{equation}
We define our cone by
\begin{equation}
  V^+(\Sol\sqc\Sol,\mathbb{Z})=
  \sum_{\mathfrak{c}\in\D(\Sol\sqc\Sol)}\mathbb{Z}^+\cdot \mathfrak{c}.
\end{equation}
A problem being given by an initial position $I$ and a final one $J$, the
\emph{simple quadratic test} consists is saying that $\1_I\sqc\1_I-\1_J\sqc\1_J\in
V^+(\Sol\sqc\Sol,\mathbb{Z})$, which can again be solved with integer
linear programming. However the spaces are much larger, and the resolution
becomes more troublesome. Note the following Lemma:
\begin{lem}
\begin{equation*}
  |\Sol\sqc\Sol|=|\Sol|(|\Sol|+1)/2\quad,\quad
  |\D(\Sol\sqc\Sol)|=(|\Sol|-2)|\D(\Sol)|.
\end{equation*}
\end{lem}
Indeed, there are $|\D(\Sol)|$ moves of type $|\cg|\sqc \cg$, and, for each
$\cg\in\D(\Sol)$, there are $|\Sol|-3$ moves of type $2\cvv{A}\sqc \cg$ with
$\cg(A)=0$. For the english board, the cardinality of $ |\D(\Sol\sqc\Sol)|$ is
thus $2\,356$ for a board of 561 squares.

We have already given an example showing that this test is sometimes
better than the linear test with non-negative integer coefficients
but we show now that this is always the case. To do so, let us define
\begin{equation*}
 \left\{
    \begin{array}{l}
      \displaystyle\F_0(\Sol\sqc\Sol,\mathbb{Z})=
      \sum_{A\in\Sol}\mathbb{Z}\cdot \cvv{A\sqc A}
      +\sum_{A\neq B\in\Sol}\mathbb{Z}\cdot 2\cvv{A\sqc B}
      \\
      \displaystyle W(\Sol\sqc\Sol,\mathbb{Z})=
      \sum_{A\neq B\in\Sol}\mathbb{Z}\cdot 2\cvv{A\sqc B}.
    \end{array}
    \right.
\end{equation*}
Then we can easily identify
$\F_0(\Sol\sqc\Sol,\mathbb{Z})/W(\Sol\sqc\Sol,\mathbb{Z})$ with the space of
integer valued functions on $\{\cvv{A\sqc A}, A\in\Sol\}$, which we can in turn
identify with $\Sol$. By these identifications, we start with a function $h\in
\F(\Sol,\mathbb{Z})$, build $h\sqc h\in \F_0(\Sol\sqc\Sol,\mathbb{Z})$ and is
next send to $h$. In particular, we get
\begin{equation}
  \label{eq:11}
  \1_I\sqc\1_I-\1_J\sqc\1_J\in V^+(\Sol\sqc\Sol,\mathbb{Z})
  \implies
  \1_I-\1_J\in V^+(\Sol,\mathbb{Z}).
\end{equation}
The fact that this test is in fact strictly superior on some boards in shown
by the problem described by figures~\ref{cvvd1} and~\ref{cvvf1}.

\section{A quadratic test, with flatness constraints}

If the simple quadratic test is stronger than the linear one with
positive integers, it turns out when used to be lacking in
efficiency. The last term in~\eqref{eq:10} can be
written as $2\1_{K}\sqc \cg$ where $K\subset\Sol$ avoids the support of
$\cg$. This is much better than saying that it is a linear combination of
$2\cvv{A}\sqc \cg$, but it leads to $2^{|\Sol|-3}|\D(\Sol)|+|\D(\Sol)|$
generators! This is of course way too much and makes this new set of
generators impractical. However, if $\Fcal$ is a succession of legal moves from $I$
to $J$, we can write
\begin{equation}
  \label{eq:12}
   \1_I\sqc\1_I-\1_J\sqc\1_J
  =
  \sum_{\cg}x(\cg) |\cg|\sqc\cg+
  \sum_{\cg}\sum_{\substack{A}}y_{\cg}(A)2\cvv{A}\sqc \cg.
\end{equation}
And we readily see that on this writing that the following
inequalities are satisfied
\begin{equation}
  \label{eq:13}
  0\le y_{\cg}(A)\le x(\cg).
\end{equation}
We call them the \emph{flatness constraints}. Despite their number, these constraints renders the
quadratic test much more efficient. In fact, The $x(\cg)$ are related to the
usual linear moves by
\begin{equation}
  \1_I-\1_J=\sum_{\cg}x(\cg) \cg
\end{equation}
(see the process that enabled us to prove~\eqref{eq:11}) 
and as such can be controlled in size by the thickness of $\cg$ at
$\1_I-\1_J$, as defined in section~\ref {thickness}.

On an english board, the $x(\cg)$'s are seldomly
larger than $4$, and on arbitrary board they are anyway bounded.

Notice that if $\1_I\sqc\1_I-\1_J\sqc\1_J$ passes this test, then actually, it can be
written as a linear combination with non-negative integer coefficients of
$2\cvv{A}\sqc \cg$ with $\cg(A)=0$ and diagonal moves $|\cg|\sqc\cg$. To
realize such a writing, given $\cg$, simply collect together all $A$'s for which
$y_{\cg}(A)$ has a given value into a set $\mathcal{A}$. Note that these
sets $\mathcal{A}$ are \emph{not} the same as the sets $K$ we used at the very
beginning of this section, but are of same use.

The problem described by figures~\ref{ex7} and~\ref{ex8}  goes through the quadratic test with no flatness
constraints, but is shown impossible as soon as we add these constraints~:

\begin{figure}[!h]
  \centering
  \begin{minipage}[c]{0.4\textwidth}
     \centering \scalebox{0.3}{\includegraphics{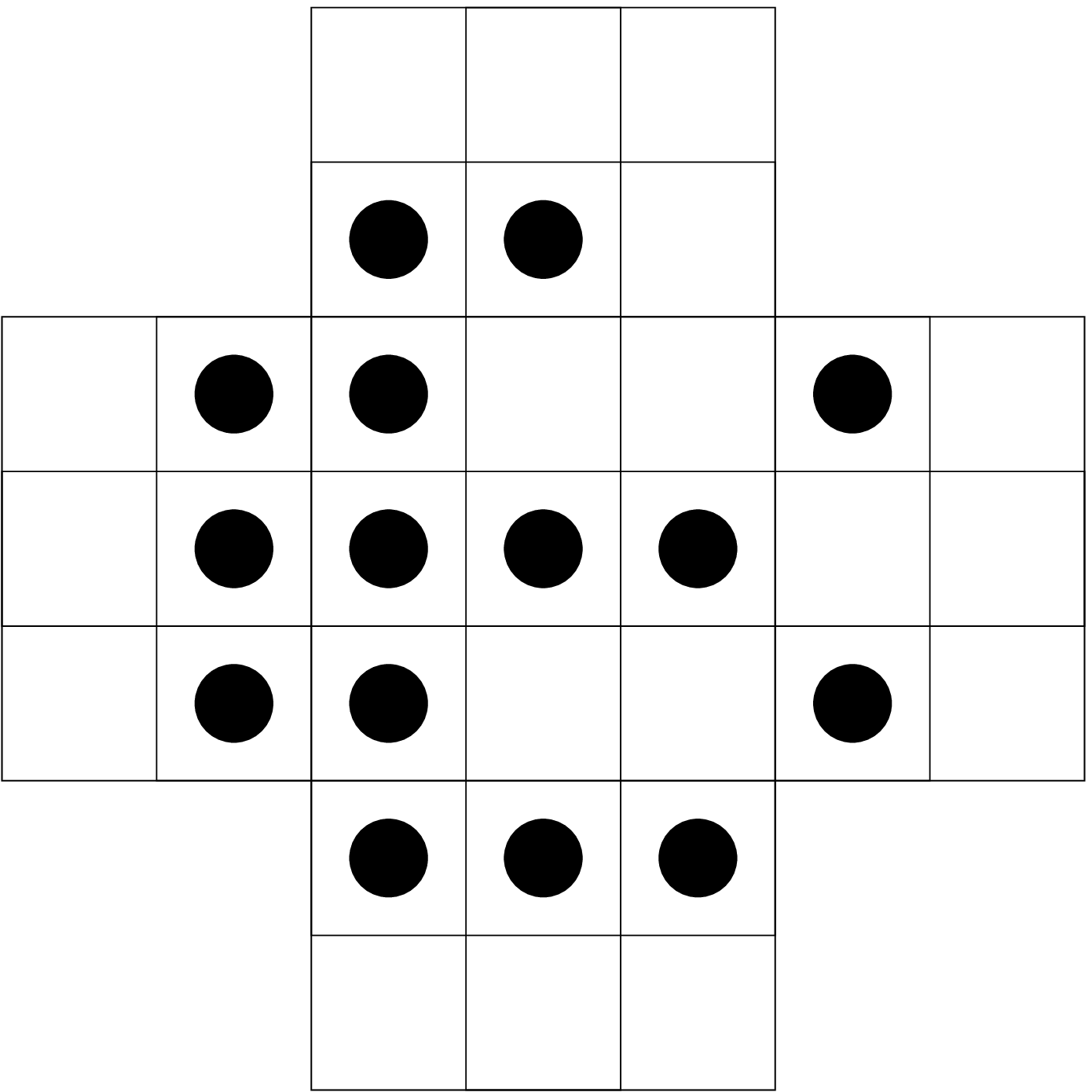}}%
     \caption{Starting position}\label{ex7}
  \end{minipage}\quad\quad
  \begin{minipage}[c]{0.4\textwidth}
     \centering \scalebox{0.3}{\includegraphics{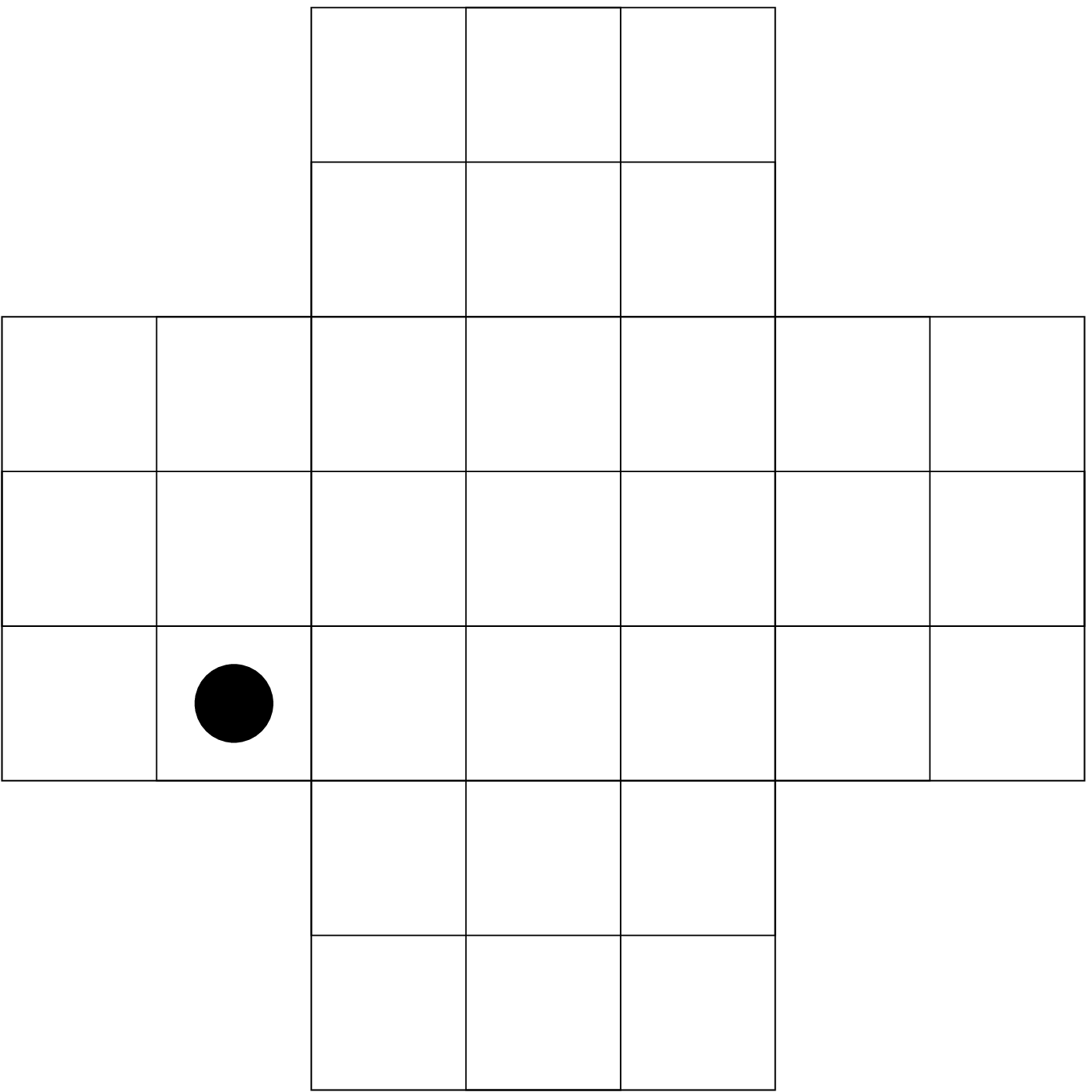}}%
     \caption{Ending position}\label{ex8}
  \end{minipage}
\end{figure}

\FloatBarrier

This new test is the main novelty of this paper and is extremely efficient in
practice, though it requires a processor to carry out the required computations.

We end this part with three further examples of problems shown to be impossible
via the quadratic test with flatness constraints. Here are two problems, with
a same starting position but different ending positions. None of them go
through the quadratic test with flatness constraints:

\begin{figure}[!h]
  \centering
  \begin{minipage}[l]{0.26\textwidth}
  \centering
  \scalebox{0.26}{\includegraphics{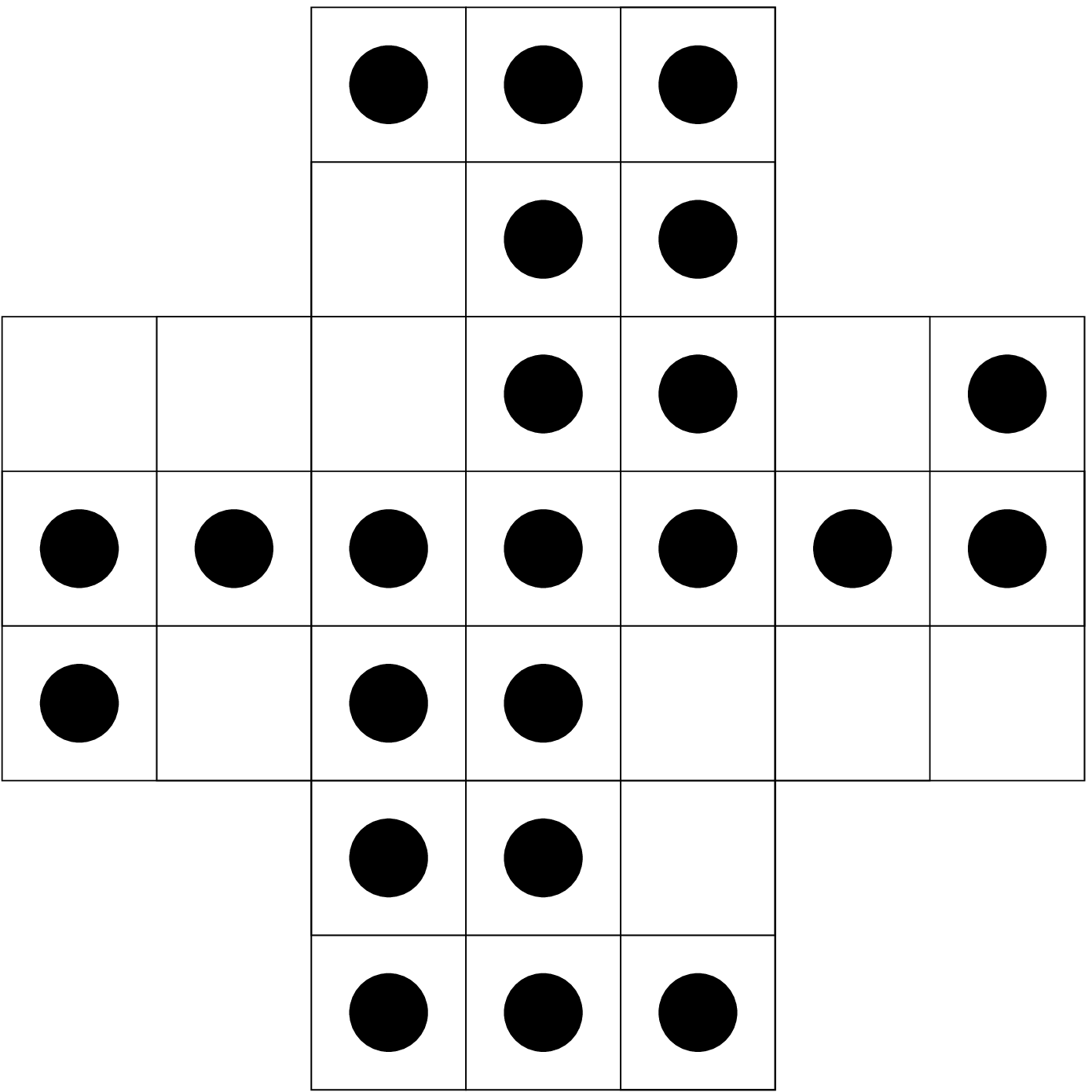}}
     \caption{Starting position}
  \end{minipage}\quad\quad
  \begin{minipage}[c]{0.26\textwidth}
  \scalebox{0.26}{\includegraphics{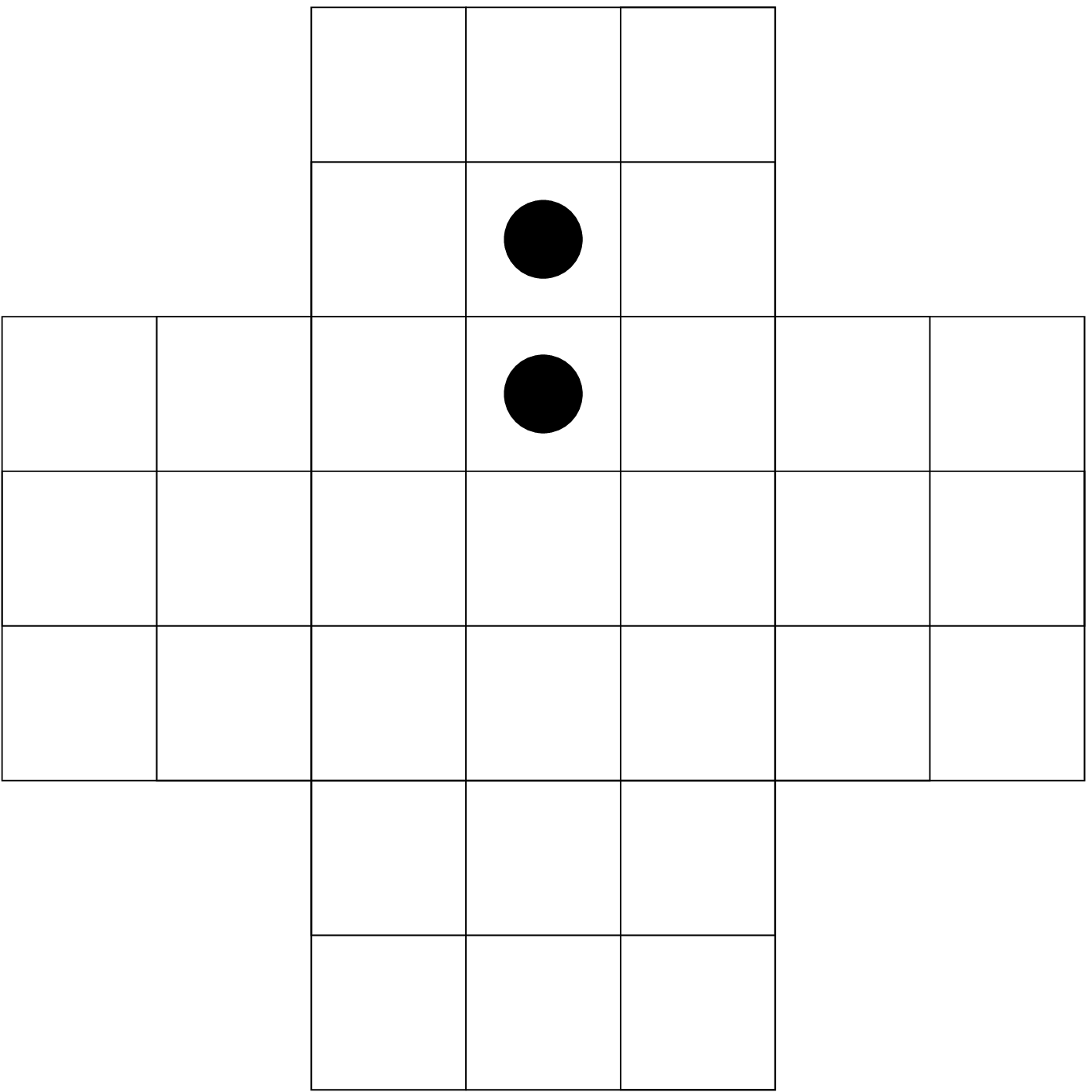}}
     \caption{First ending position}
  \end{minipage}\quad\quad
  \begin{minipage}[r]{0.26\textwidth}
  \scalebox{0.26}{\includegraphics{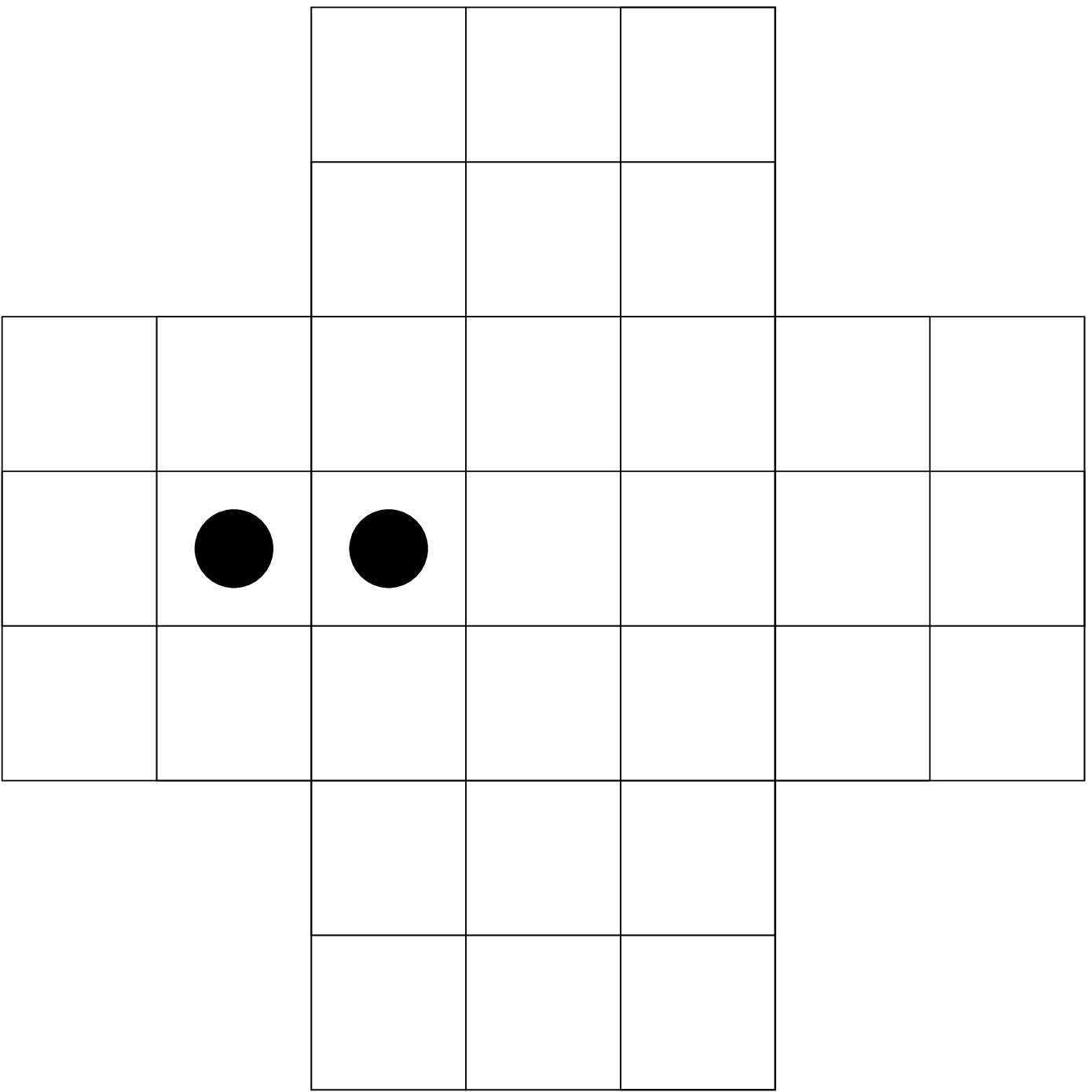}}
     \caption{Second ending position}
  \end{minipage}
\end{figure}
\FloatBarrier

The third example is to go from the initial position to the intermediate
ending position. This is shown to be impossible via the quadratic test with
flatness constraints, though it again passes the simple quadratic
test. Moreover, the problem to go from the initial position to the final
ending position is feasible in legal moves.

\begin{figure}[!htbp]
  \centering
  \begin{minipage}[c]{0.26\textwidth}
  \centering
  \scalebox{0.26}{\includegraphics{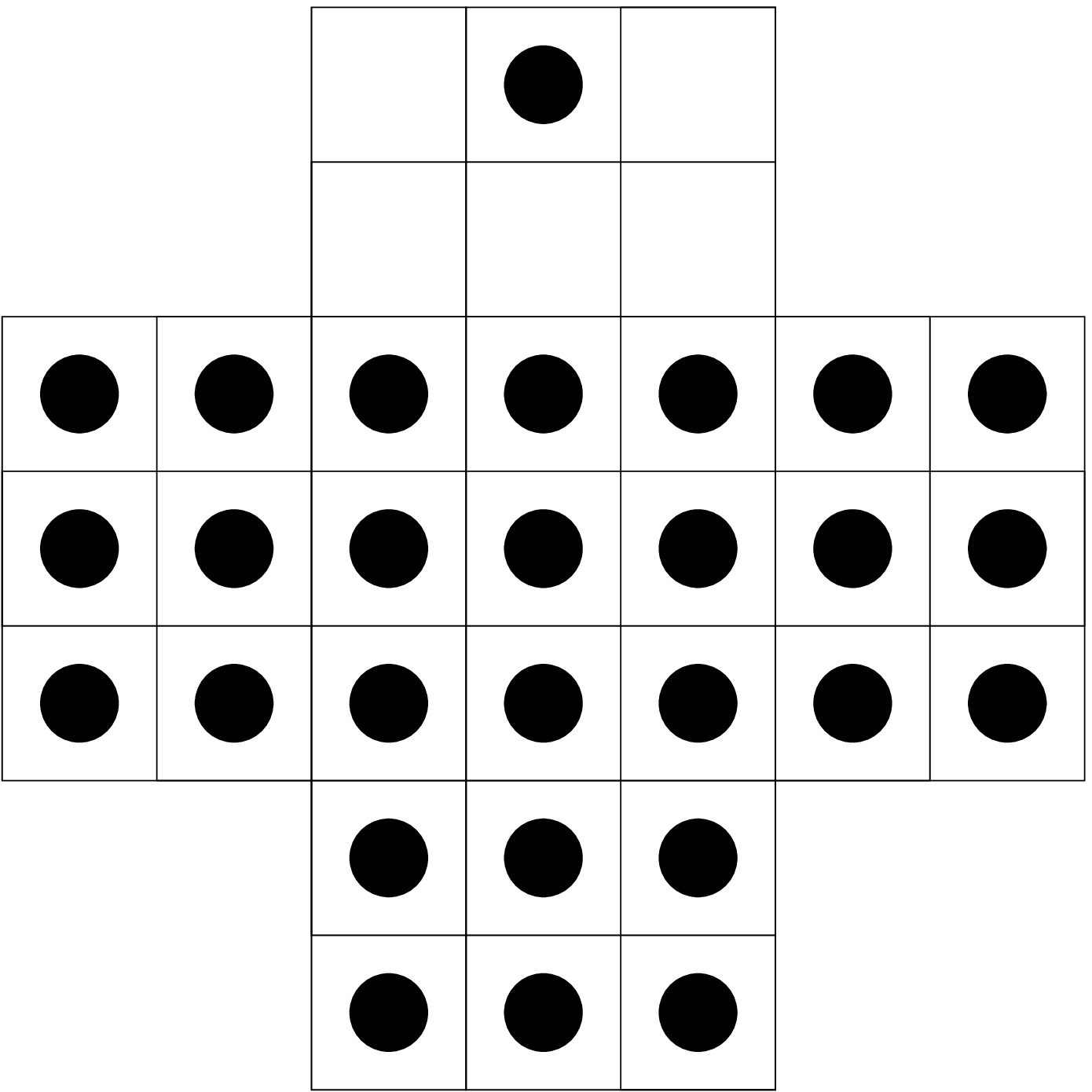}}
     \caption{Starting position}
  \end{minipage}\quad\quad
  \begin{minipage}[c]{0.26\textwidth}
   \scalebox{0.26}{\includegraphics{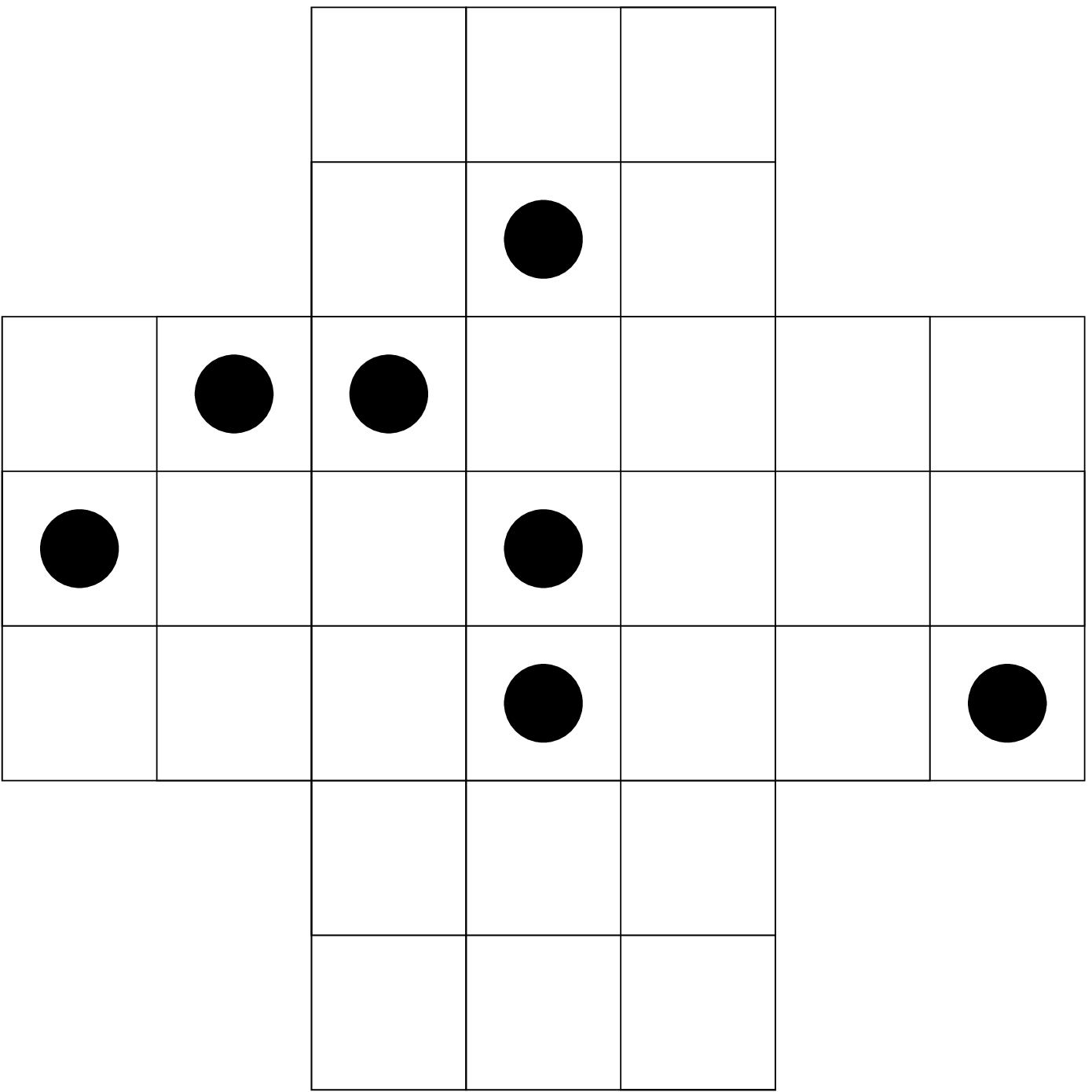}}
    \caption{Intermediate ending position}
  \end{minipage}\quad\quad
  \begin{minipage}[c]{0.26\textwidth}
  \scalebox{0.26}{\includegraphics{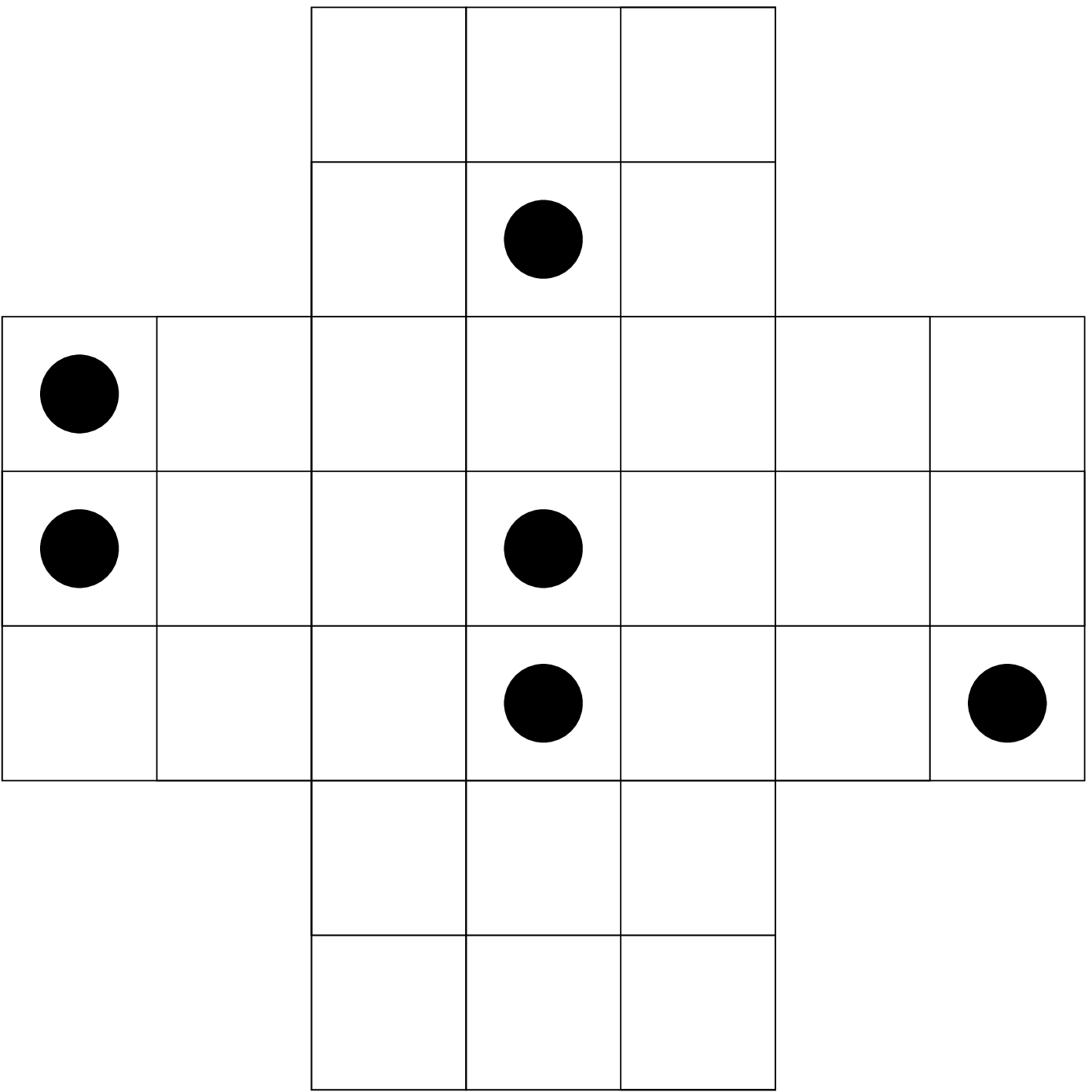}}
     \caption{Final ending position}
  \end{minipage}
\end{figure}

\FloatBarrier

\section{Additional constraints, a first draft}

Now that we have seen that the quadratic test with flatness constraints is so
very efficient, it is tempting to try to add some further constraints. This is
the topic of these two last sections, but this part is still very much in
progress. The reader may get the impression that it is not so much in
progress than more bluntly unfinished. After some months of efforts, I have not been
able to derive a unifying setup for what look like protrusions of a
hidden structure, which is why I deliver them in tnhis state.

The idea we follow is to add geometrical information to control as much as
possible these new variables $y_\cg(A)$ in \eqref{eq:12}.

Let us start with a fundamental inequality.
\begin{prop}
  Assume we can go from $I$ to $J$ in legal moves. Then there exists a writing
  of $\1_I\sqc\1_I-\1_J\sqc\1_J$ (as in~\eqref{eq:12}) such that for every
  $A\in\Sol$ we have
  \begin{equation}\label{fund}
      0\le \sum_{\cg}y_{\cg}(A)+\sum_{\cg/\cg(A)\neq0}x(\cg)\le |I|-|J|.
  \end{equation}
\end{prop}
See~\eqref{eqz2} and \eqref{eq:18} for refinements.
Let $\Fcal$ be a succession of legal moves from $I$ to $J$. We set
\begin{equation}
  \label{eq:14}
  \pfrak(A, \Fcal)=\sum_{\cg}y_{\cg}(A)+\sum_{\cg/\cg(A)\neq0}x(\cg)
\end{equation}
where the $y_{\cg}(A)$'s and the $x(\cg)$'s come from  \eqref{eq:12}.
\begin{dem}
  Given a move $\cg$, let us look at the situation of the board before using
  this move. There are four possibilities for $A$:
  \begin{itemize}
  \item $\cg(A)=1$, which means that $A$ is on the board and
    participates to the move. It is counted in $x(\cg)$ and nowhere else.
  \item $A$ is not on the board but is created by the move. It is
    counted in $x(\cg)$ and nowhere else.
  \item $A$ is on the board but does not participate to the move. It
    is counted in $y_\cg(A)$ and nowhere else.
  \item $A$ is not on the board and not created by the move. It is not
    counted anywhere.
  \end{itemize}
  The proof
  follows by using this remark and an induction on $|I|-|J|$.
  We have equality if and only if the last case above never occurs,
  which means that $A$ is never absent from the position
  for two consecutive moves.
\end{dem}

We have seen that we can have equality in~\eqref{fund}, but we can even show
that the right hand side is on average of the correct order of
magnitude. Indeed we have
\begin{eqnarray*}
  \sum_{A\in\Sol}
   \pfrak(A,\Fcal)
  &=&|I|-2+|I|-3+\cdots+|J|-1+3(|I|-|J|)
  \\&=&(|I|-|J|)\frac{|I|+|J|+3}{2}
\end{eqnarray*}
since there are $|I|-2$ points on the first move that are on the board but do
not participate to the move, then $|I|-3$, and so on. As a consequence
\begin{equation*}
  \frac{1}{|\Sol|}\sum_{A\in\Sol}
  \pfrak(A,\Fcal)
  =\bigl(|I|-|J|\bigr)\frac{|I|+|J|+3}{2|\Sol|}.
\end{equation*}
This shows that~\eqref{fund} prevents too wide deviations from the mean, at
least if $|I|+|J|$ and $|\Sol|$ are of comparable size. We propose to improve
on this double inequality in three ways.

\subsection{Using the speed at which a peg gets inside $J$}
We define the \emph{depth} of the
point~$A$ with respect to the position $S$ containing it to be the minimum
number $\Depth(A,S)$ of legal moves
required to remove the peg in~$A$. If $A$ is not in $S$, we set
$\Depth(A,S)=0$. Let us recall a classical Lemma.
\begin{lem}[Leibniz]\label{Leibniz}
  If the sequence of legal moves
  $\cg_1,\cg_2,\dots,\cg_k$ goes from $I$ to $J$, then the sequence of legal
  moves $\cg_k,\dots,\cg_2,\cg_1$ goes from $\Sol\setminus J$ to
  $\Sol\setminus I$.
\end{lem}
It is enough to verify this property when $k=1$ where it is
obvious. Leibniz expressed this idea in a different manner: he started from
the final position $J$ and tried to recover the initial one by playing in
reverse; he discovered it was the same game, provided one considered the empty
squares as having a peg, and the ones with a peg as being empty. This is
exactly what we shall consider. Indeed, given a point $A$ out of our final
position~$J$, there is a minimal number a moves that will "bring" its peg
inside $J$, or kill it, namely $\Depth(A,\Sol\setminus J)$.

Let us select a minimal path from $J$ to $A$. Its last move puts a peg in $A$, 
i.e. has $A$ as point $R$ since we could otherwise shorten this path. Moreover
it does not use $A$ anymore as point $P$ or $Q$ since we 
could again shorten the path. Consequently, for any $A\notin J$
\begin{equation}
  \label{eq:15}
  \pfrak(A,\Fcal)\le \max(0,|I|-|J|-\Depth(A,\Sol\setminus J)+1).
\end{equation}
If $A$ is in $J$, we have $\Depth(A,\Sol\setminus J)=0$ so that~\eqref{fund}
is stronger.
\begin{dem}
  Indeed  $A$  not in $J$ implies $\Depth(A,\Sol\setminus J)\ge1$.
  The $\Depth(A,\Sol\setminus J)-1$ last moves cannot use $A$ in any part of a move, hence
  we can use~\eqref{fund} with $|J|+\Depth(A,\Sol\setminus J)-1$ points as a final
  position instead of $J$ if $A$ is at some point of time on the board. Else,
  it is never here and the upper bound 0 is fine.
\end{dem}
We do not know of any precise mean of computing this depth, but we provide
now a fast way to get an excellent lower bound. Let us consider the
oriented graph $\mathfrak{G}$ built on the set $\Sol$ and where we put an edge
from $A$ to $B$ if there exists $\cg\in\D(\Sol)$ such that $\cg(A)=1$ and
$\cg(B)=-1$. A minimal path that realizes $\Depth(A,\Sol\setminus J)$ is readily
transformed in a path from $A$ to $J$ on $\mathfrak{G}$. Reciproquely
from such a path from $A$ to $J$ on this graph, we deduce a position
$K$ by adding the required points $P$ and $Q$ necessary
for the $\cg$'s. The only problem is that this process may require to put
several pegs on a same square (we do not have any example of such a situation).
Denoting by
$\delta_{\mathfrak{G}}(A,J)$ the distance on this graph, we have established that
\begin{equation}
  \label{eq:16}
  \delta_{\mathfrak{G}}(A,J)\le \Depth(\Sol\setminus J)
\end{equation}
Note that a final position $L$ in case of $\delta_{\mathfrak{G}}$ is reduced
to a single point. The distance $ \delta_{\mathfrak{G}}(A,J)$ is now readily
computed, by using the Dijkstra's algorithm for instance.

Practically, to find a minorant of this depth, we proceed in two steps (with
$S=\Sol\setminus J$):
\begin{itemize}
\item We try every succession of 5 legal moves from $S$.
\item Concerning the remaining ones, we first build the set $S_5$ of points with
  $\Depth(A,S)\le 5$. If $A\in S_5$, we find the minimum of
  $\delta_{\mathfrak{G}}(A,B)+\Depth(B,S)$ for every $B\in S_5$; this a first
  lower bound for $\Depth(A,S)$, but sometimes the lower bound~6 is simply better.
\end{itemize}

\subsection{Using the speed at which a point is reached by $I$}
Let us now examine the somewhat reciproqual situation, and try to get the 
minimum of legal moves from the set $I$ that puts a peg in $A$.
We need two pegs to create one, which means that the distance
$\delta_{\mathfrak{G}}(A,I)$ is not a good lower bound anymore.
We define the \emph{height} $\Height(A,I)$ of A with respect to $I$ to be his
minimal number, and set $\Height(A,I)=\infty$ if $A$ can never be reached. 
Computing this $\Height(A,I)$ is very difficult.

\begin{lem}\label{Horsconvexe}
  Let $I$ be a subset of $\Sol$ and $A$ be such that
  $\Height(A,I)<\infty$. For any non-negative resource count
  $\pi$, we have
  $\langle{\1_I},{\pi}\rangle\ge \pi(A)$.
\end{lem}

\begin{dem}
  Indeed there is a set $J$ which contains $A$ and that is reachable from
  $I$. We thus have $\langle{\1_I},{\pi}\rangle\ge \langle{\1_J},{\pi}\rangle$
  which in turn is non less than $\pi(P)$ by the non-negativity assumption
  on $\pi$.
\end{dem}

Using Lemma~\ref{Horsconvexe} and some direct computations, we get the
following height-diagram for the left-hand side position. 

\begin{figure}[!h]
  \centering
  \begin{minipage}[l]{0.32\textwidth}
  \centering
  \scalebox{0.32}{\includegraphics{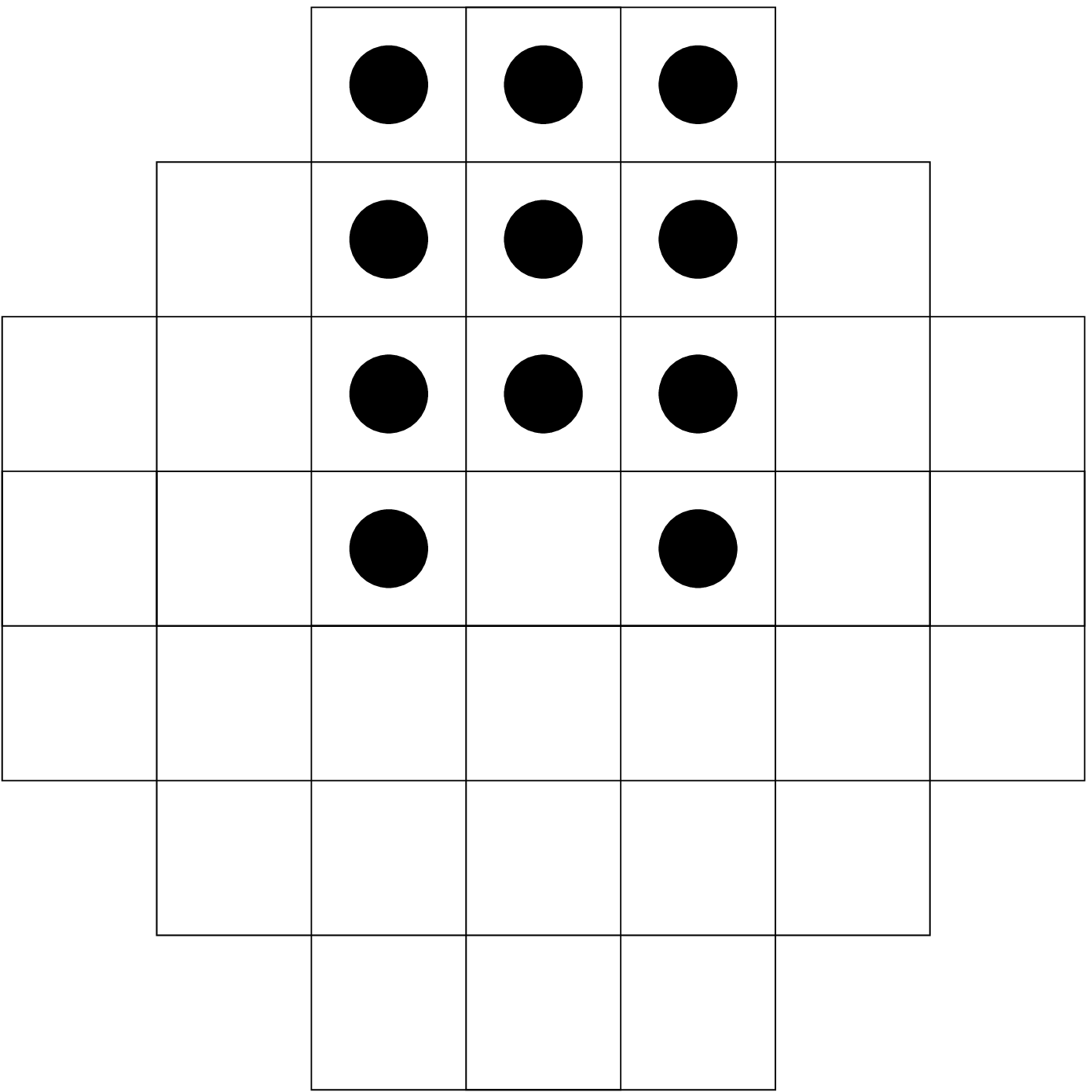}}
     \caption{Starting position}
  \end{minipage}\quad\quad
  \begin{minipage}[r]{0.32\textwidth}
  \scalebox{0.32}{\includegraphics{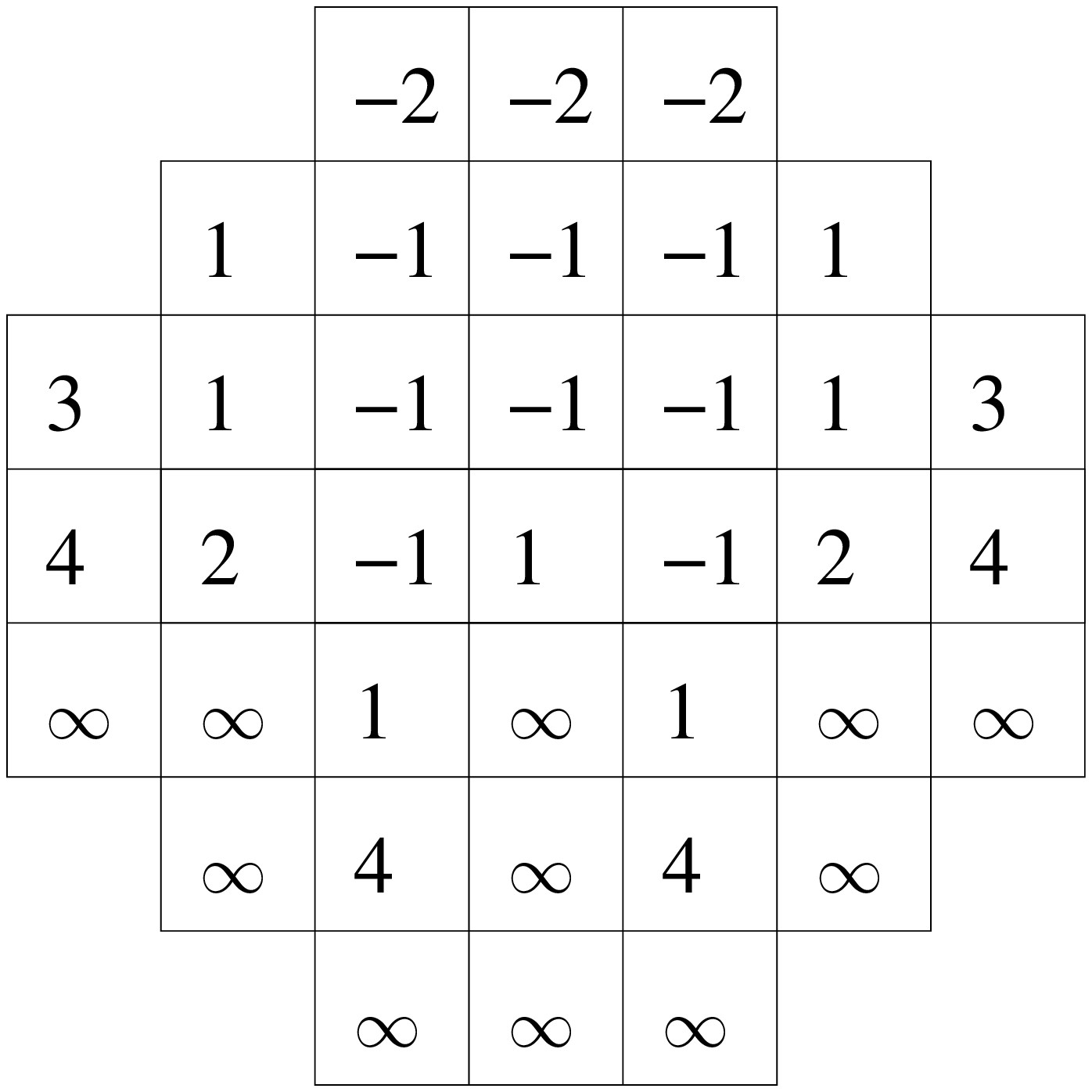}}
     \caption{Height / -Depth }
  \end{minipage}
\end{figure}
\FloatBarrier

We next provide an example on which Lemma~\ref{Horsconvexe} is not strong
enough to decide whether some points have finite heights or not.
This problem passes the linear integer test. We provide the height of
each square (we simply computed all position attainable in 5 moves~!). The two
squares on the left-hand side (and the symmetric ones on the right-hand side)
are rather clearly not reachable, but the test deduced from
Lemma~\ref{Horsconvexe} fails to prove that. Even worse, we found for each of
this square a position got from the first one in 5 moves and for which this
square is not shown to be unreachable by this test.

\begin{figure}[!h]
  \centering
  \begin{minipage}[l]{0.26\textwidth}
  \centering
  \scalebox{0.20}{\includegraphics{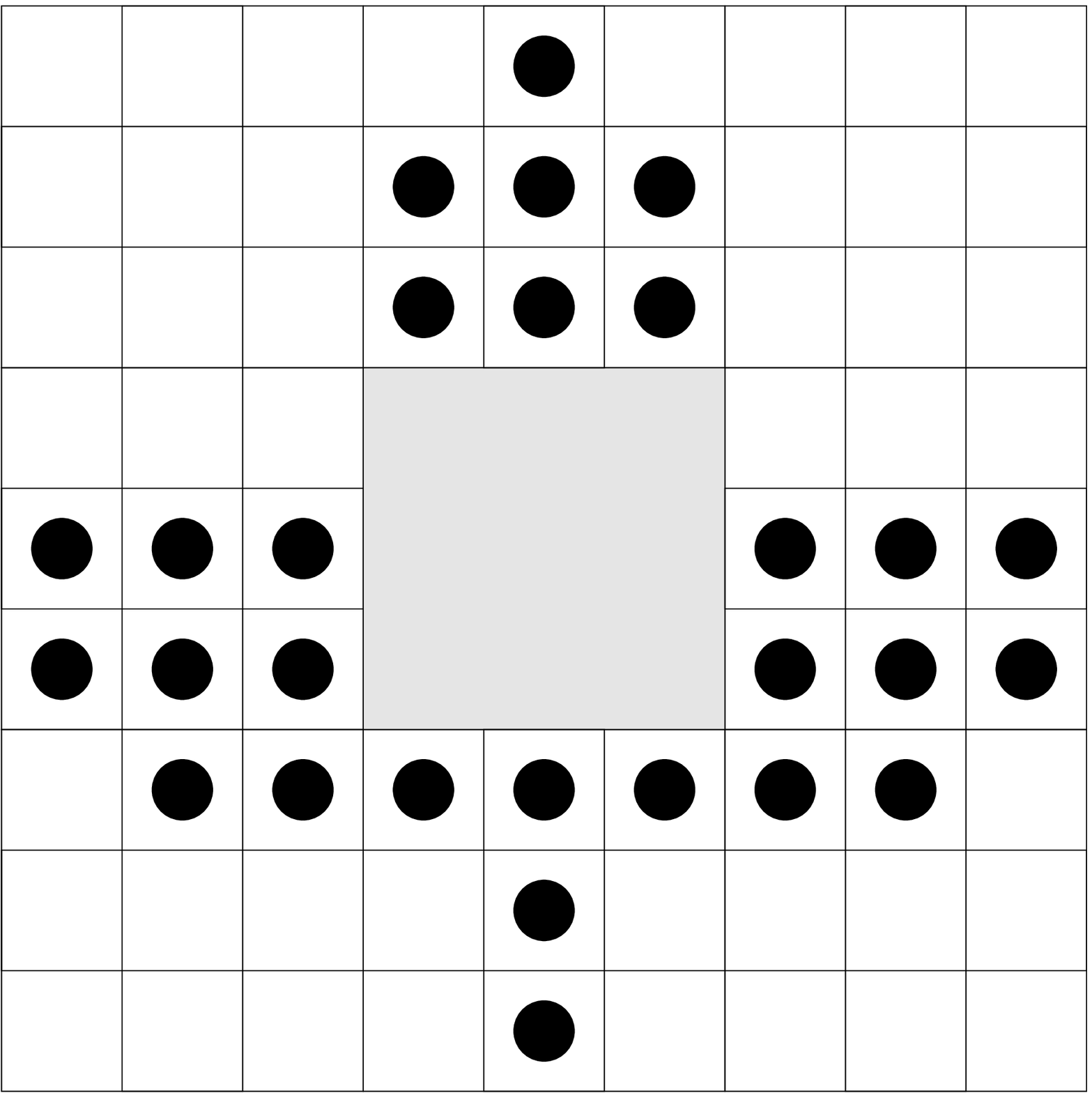}}
     \caption{Starting position}\label{ce3}
  \end{minipage}\quad\quad
  \begin{minipage}[c]{0.26\textwidth}
  \scalebox{0.20}{\includegraphics{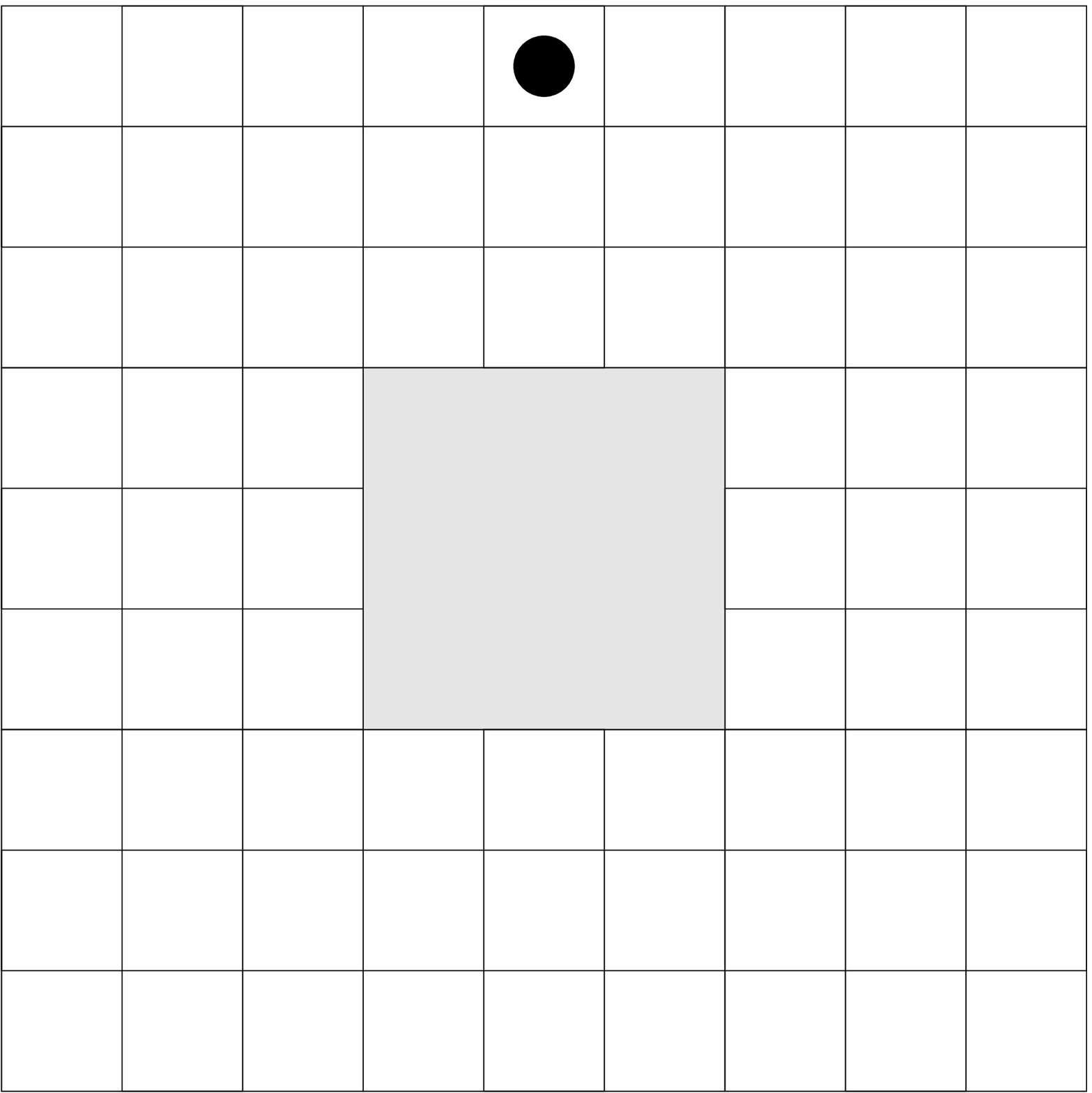}}
     \caption{Ending position}\label{ce4}
  \end{minipage}\quad\quad
  \begin{minipage}[r]{0.26\textwidth}
  \scalebox{0.20}{\includegraphics{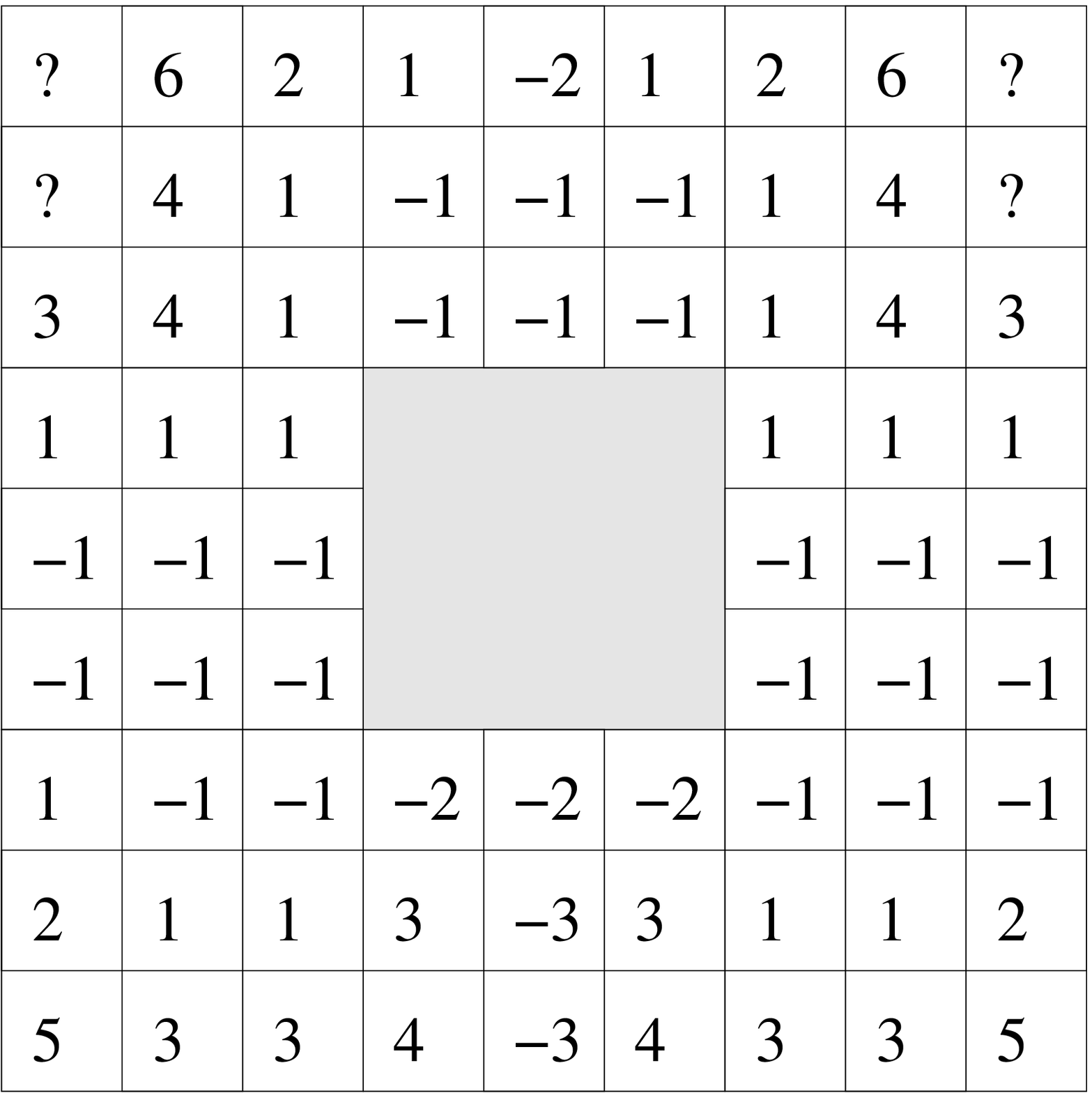}}
     \caption{Height / -Depth}
  \end{minipage}
\end{figure}
\FloatBarrier

Practically, to find a minorant of this height, we proceed in three steps:
\begin{itemize}
\item We try every succession of 5 legal moves.
\item We use the Lemma~\ref{Horsconvexe} to determine those points that are
  guaranteed to have infinite height. (We apply this test to all of the
  derived positions).
\item Concerning the remaining ones, we first build the set $I_5$ of points with
  $\Height(A,I)\le 5$. If $A\in\Sol\setminus I_5$, we find the minimum of
  $\delta_{\mathfrak{G}}(A,B)+\Height(B,I)$ for every $B\in I_5$; this a first
  lower bound for $\Height(A,I)$, but sometimes the lower bound~6 is simply better.
\end{itemize}

Set
\begin{multline}
  \label{eq:17}
  \corr(A, I, J)=\min \Bigl(|I|-|J|,
  \max\bigl(\Depth(A,\Sol\setminus J)-1, 0\bigr)
  \\+\max\bigl(\Height(A,I)-1, 0\bigr)\Bigr).
\end{multline}
We have
\begin{equation}
  \label{eqz2}
  \pfrak(A,\Fcal)\le |I|-|J|-\corr(A, I, J).
\end{equation}

\subsection{Using the speed at which a peg comes out of $I$}
We finally improve on the lower bound in~\eqref{fund}. The fact is that some
points are so much within the starting position $I$ that the peg on them cannot
be eliminated before so many moves, and this is precisely how we defined
$\Depth(A,I)$.
We have then
\begin{equation}
  \label{eq:18}
  \Depth(A,I)\le \pfrak(A,\Fcal).
\end{equation}

\subsection{Final discussion}
We end this section with two remarks. First, both notions of depth and height
use only one of the two positions of the problem, and this is a loss. For
instance concerning height, if we manage to put a peg in a very far away
square that is also far from our final position, it is probable that we shall
not be able to bring it back to it; for instance, if the starting position is
given by figure~\ref{ce3}, it is likely that we cannot put a point in
the lower left corner and finish as in figure~\ref{ce4}. Secondly, constraints \eqref{eqz2} and
\eqref{eq:18} only avoid extremal cases, as we noted earlier, and there are
only $2|\Sol|$ of them for a problem with about $|\Sol|^2$ variables; in fact,
if $\Sol$ has no isolated point, Lemma~\ref{count} yields
\begin{equation*}
  |\Sol|(|\Sol|-2)\le |\D(\Sol\sqc\Sol)|\le 4|\Sol|(|\Sol|-2).
\end{equation*}
This explains why these constraints
are somewhat weak.

\section{Additional constraints}
Having in mind the counting argument displayed at the end of last section, we
see that finding conditions on couples $(A, A')$ of points would not increase
too much the size of the problem but may yield more stringent constraints.

As of now, we have only found one such type of constraint, which applies to
initial positions $I$ such that $\Sol\setminus I$ is large enough.

Let us start with some general considerations.  Let $\Height(A, A',I)$
be the minimum number of legal moves necessary to put a peg in each of
$A$ and $A'$, starting from a board with pegs on all the points of
$I$. We assign it value $\infty$ if no such succession exists. Note
that the height-function does not behave like a distance, since we can
have $\Height(A, A',I)>\Height(A,I)+\Height(A',I)$. We formulate a
conjecture:
\begin{conj}
  $\Height(A, A',I)\ge\Height(A,I)+\Height(A',I)$.
\end{conj}
A proof or disproof of this conjecture has sofar escaped the author.

\begin{lem}\label{addc}
  Consider two points $A$ and $A'$ such that $\Height(A, A',I)=\infty$. 
  Then
  \begin{equation}
    \pfrak(A,\Fcal)+\pfrak(A',\Fcal)-\sum_{\cg/\cg(A)\cg(A')\neq0}x(\cg)\le |I|-|J|.
  \end{equation}
\end{lem}
\begin{dem}
  Given a move $\cg$, let us look at the situation of the board before using
  this move. There are several cases:
  \begin{itemize}
  \item $A$ is on the board and is not moved by $\cg$.
    Then $A'$ is not on the board, and may not be created by $\cg$.
    This move is counted in $y_\cg(A)$.
  \item $A'$ is on the board and is not moved by $\cg$.
    Then $A$ is not on the board, and may not be created by $\cg$.
    This move is counted in $y_\cg(A')$.
  \item $A$ is on the board and is moved by $\cg$. Then $A'$ is not on the
    board and may be created. This move is counted in $x(\cg)$.
  \item $A'$ is on the board and is moved by $\cg$. Then $A$ is not on the
    board and may be created. This move is counted in $x(\cg)$.
  \end{itemize}
\end{dem}
The question arises as to whether this Lemma leads or not to improvements, and
we provide an example below showing that it indeed does. The geometrical fact
that we have used is that a square can either contain a peg, or be empty, a
fairly trivial information that was until now absent from our discussion.

Before exposing our example, let us address rapidly the problem of computing
couples $(A,A')$ with $\Height(A,A',I)=\infty$.
\begin{lem}\label{Horsconvexebis}
  Let $I$ be a subset of $\Sol$ and $A$ and $A'$ be two points of $\Sol$.
  If there exists a non-negative resource count
  $\pi$, such that
  $\langle{\1_I},{\pi}\rangle< \pi(A)+\pi(A')$, then $\Height(A, A',I)=\infty$.
\end{lem}
We can improve on this criteria: simply form all positions derived from $I$ by
one (or any fixed number) legal move, and apply this criteria to each of them.

Here is a problem that is shown impossible by using this criteria,
though it passes the quadratic integer test with flatness
constraints:
\begin{figure}[!h]
  \centering
  \begin{minipage}[l]{0.26\textwidth}
  \centering
  \scalebox{0.20}{\includegraphics{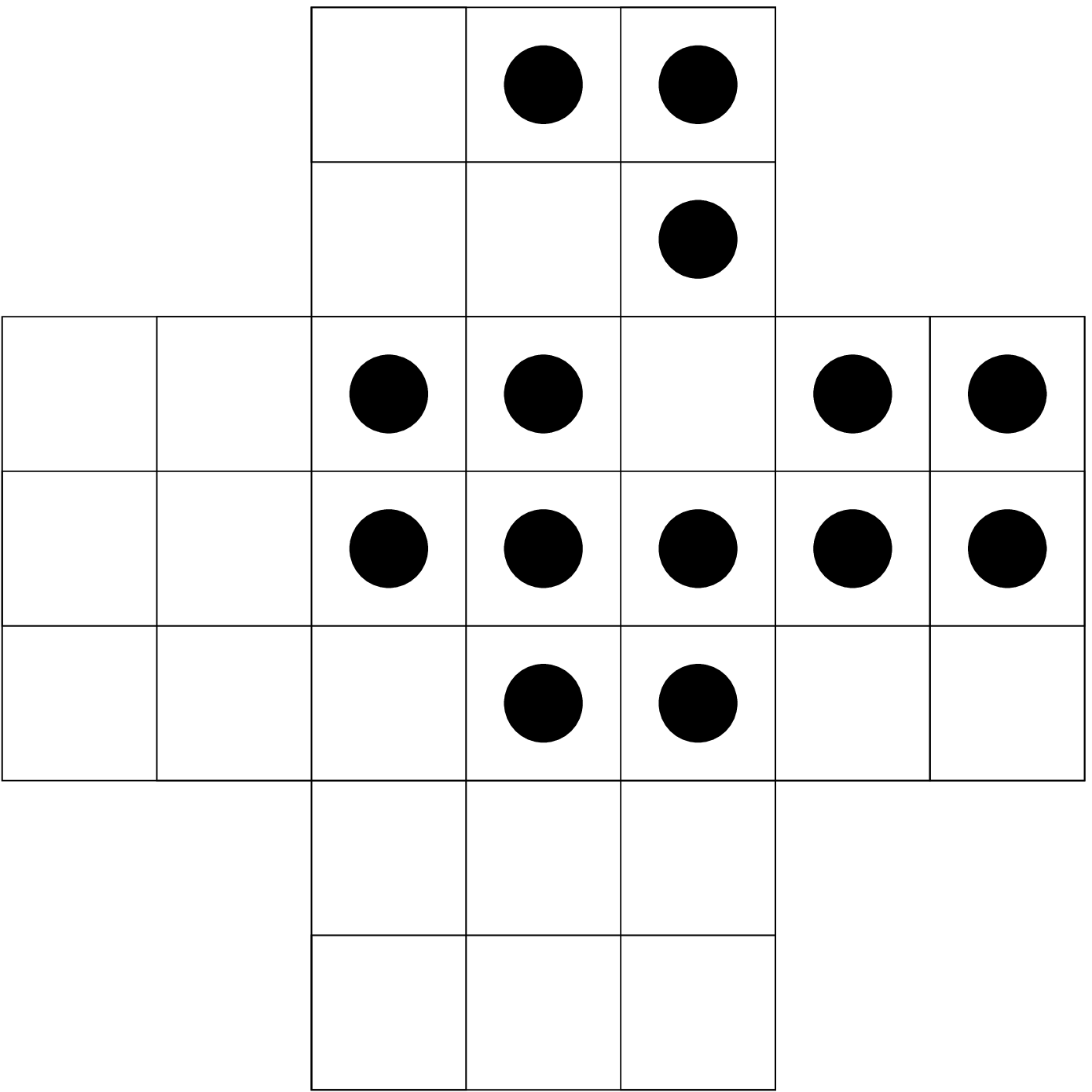}}
     \caption{Starting position}\label{ce31}
  \end{minipage}\quad\quad
  \begin{minipage}[c]{0.26\textwidth}
  \scalebox{0.20}{\includegraphics{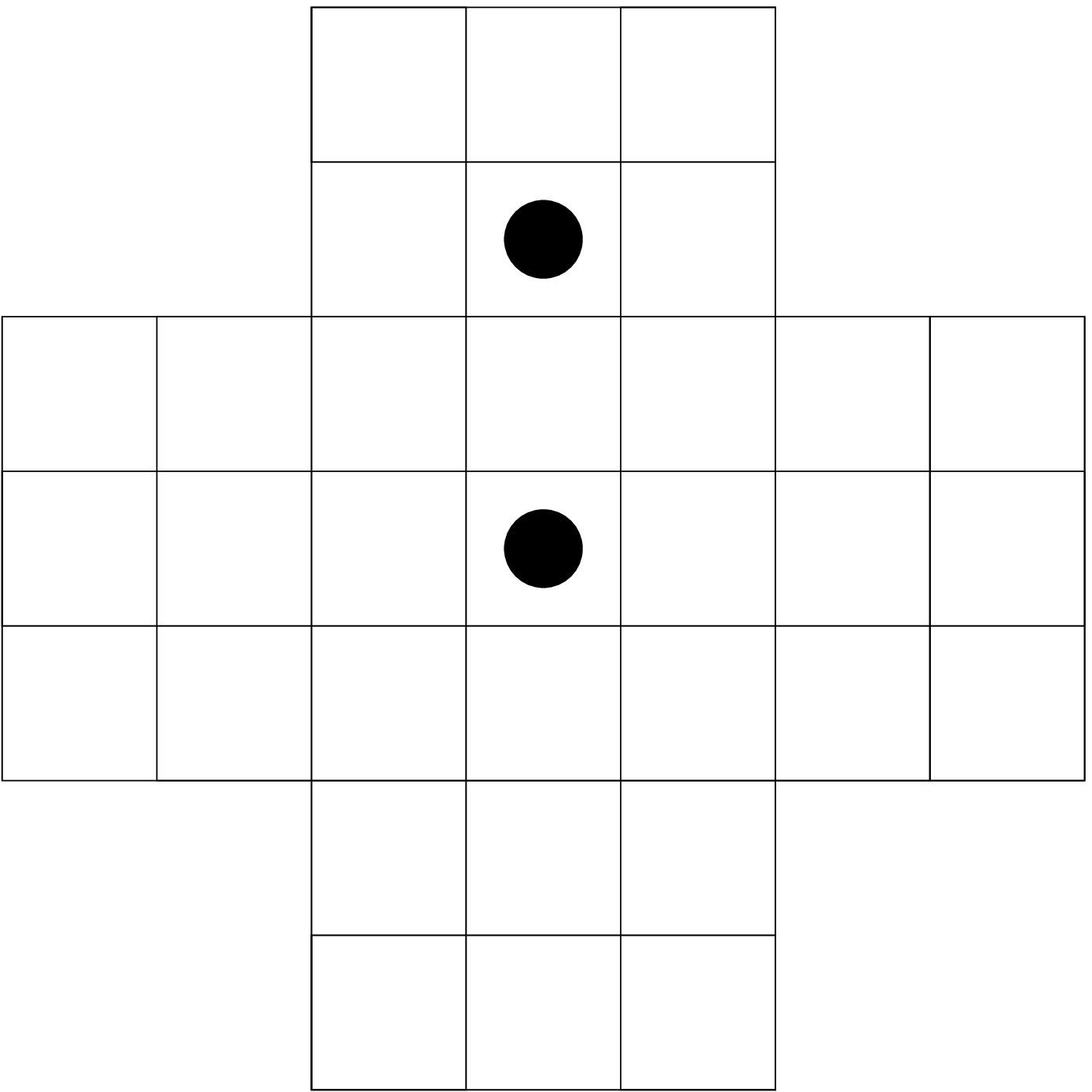}}
     \caption{Ending position}\label{ce41}
  \end{minipage}\quad\quad
  \begin{minipage}[r]{0.26\textwidth}
  \scalebox{0.20}{\includegraphics{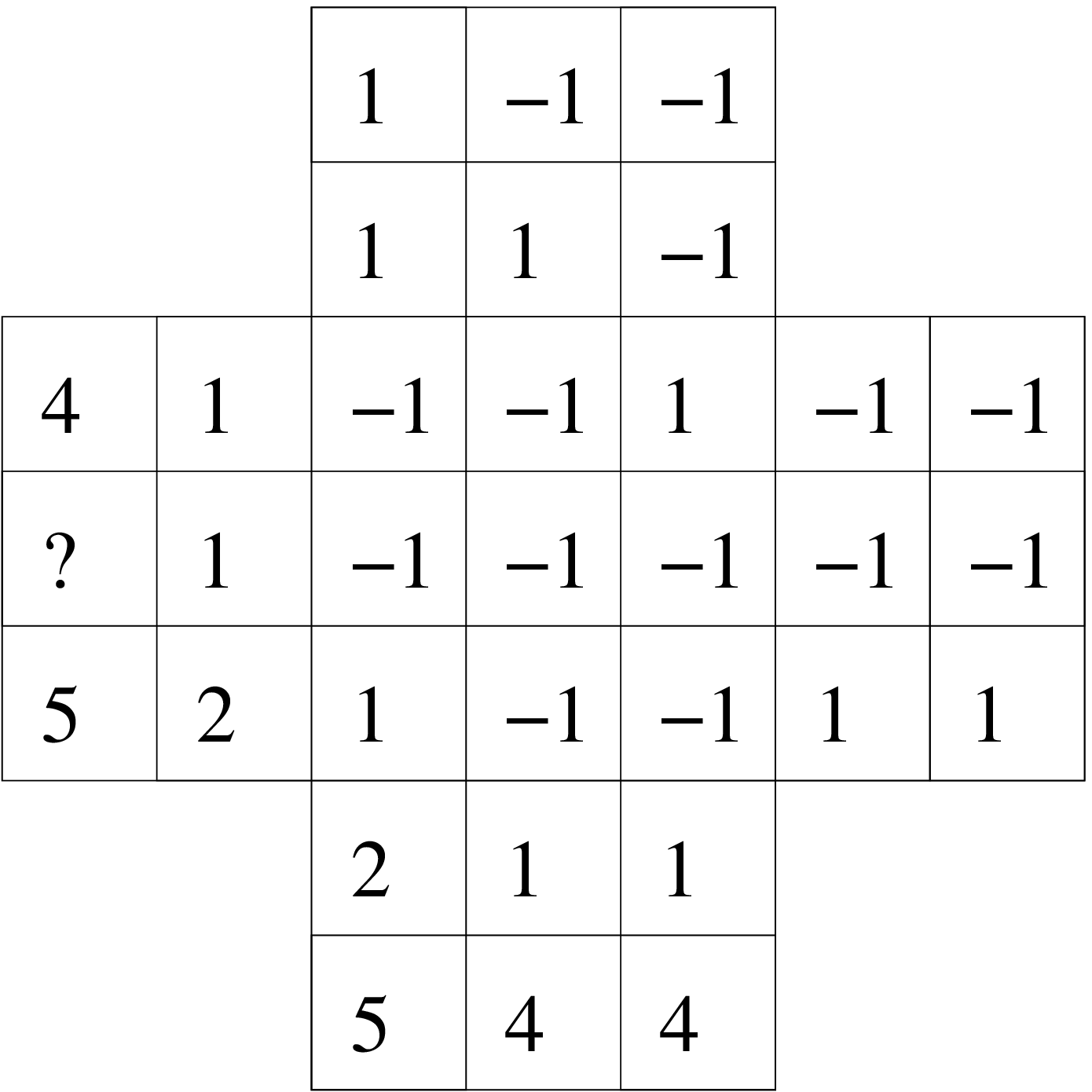}}
     \caption{Height / -Depth}
  \end{minipage}
\end{figure}
\FloatBarrier

This example is also interesting because of the square with an
interrogation dot: it is "clearly" of infinite height, but our automatic process
is not able to conclude. Here is the list of couples with
$\Height(A,A',I)=\infty$ that we have found:
\begin{equation*}
\begin{array}{|r|l|}
  \hline
  \hfill A\hfill & \hfill A'\hfill\hfill\\
  \hline
  (3, 1) & (1, 4), (2, 3), (3, 2), (3, 3)\\
  (4, 1) & (1, 4), (1, 5), (2, 3), (2, 4), (2, 5), (3, 1), (3, 3), (3, 4), (3, 6), 
  \\&(3, 7), (4, 2), (4, 3), (4, 4), (4, 5), (4, 6), (4, 7)\\
  (5, 1) & (1, 4), (2, 3), (3, 1), (3, 2), (3, 3), (4, 1), (4, 2), (5, 2), (5, 3)\\
  (5, 2) & (3, 1), (4, 1)\\
  (5, 3) & (3, 1)\\
  (5, 4) & (4, 1)\\
  (6, 3) & (3, 1), (4, 1), (5, 1)\\
  (6, 4) & (4, 1), (5, 1)\\
  (6, 5) & (4, 1)\\
  (7, 3) & (1, 4), (2, 3), (3, 1), (3, 2), (4, 1), (4, 2), (4, 3), (4, 6), (5, 1), 
  \\&(5, 2), (5, 3), (6, 3), (6, 4), (6, 5), (7, 4), (7, 5)\\
  (7, 4) & (1, 4), (2, 4), (3, 1), (4, 1), (4, 2), (4, 4), (4, 6), (4, 7), (5, 1),
  \\& (5, 2), (5, 4), (5, 6), (6, 3), (6, 4), (7, 5)\\
  (7, 5) & (4, 1), (5, 1), (6, 3), (6, 5)
  \\\hline
\end{array}
\end{equation*}

Lemma~\ref{addc} is of course of fairly limited use: we need the starting
position to leave free enough squares on the board. However, it shows how more
geometrical arguments may be used to get improvements! Our journey ends here.


\begin{thebibliography}{10}
\expandafter\ifx\csname url\endcsname\relax
  \def\url#1{\texttt{#1}}\fi
\expandafter\ifx\csname urlprefix\endcsname\relax\def\urlprefix{URL }\fi

\bibitem{Leibniz*10}
G.~Leibniz, Annotatio de quibusdam ludis, M\'emoire de l'Acad\'emie des
  Sciences de Berlin (Miscellane Berolensia).

\bibitem{Beasley*92}
J.~Beasley, Ins and Outs of Peg Solitaire, Recreations in Mathematics Series,
  Oxford University Press, 1985, (paperback Edition 1992, contain an additional
  page: Recent Developments).

\bibitem{Harang*97}
E.~Harang, \url{http://eternitygames.free.fr/Solitaire_english.html} (1997).

\bibitem{DosSantos*99}
F.~{Dos Santos}, \url{http://dauphinelle.free.fr/solitaire/} (1999).

\bibitem{Uehara-Iwata*90}
R.~Uehara, S.~Iwata, Generalized {Hi-Q} is {NP}-{C}omplete, Trans IEICE 73
  (1990) 270--273.

\bibitem{Ravikumar*97}
B.~Ravikumar, Peg-solitaire, {S}tring {R}ewriting {S}ystems and {F}inite
  {A}utomata, in: Springer (Ed.), Proc. 8th Int. Symp. Algorithms and
  Computation, Vol. 1350 of Lecture Notes in Computer Science, 1997, pp.
  233--242.

\bibitem{Reiss*57}
M.~Reiss, Beitr{\"a}ge zur {T}heorie des {S}olit{\"a}r-{S}piels, Crelles
  Journal 54 (1857) 344--379.

\bibitem{Lucas*91}
E.~Lucas, {R\'ecr\'eations} {Math\'ematiques}, 2nd Edition, Gauthiers Villars
  et fils, imprimeurs-libraires, 1891, 87-141.

\bibitem{Berlekamp-Conway-Guy*82}
E.~Berlekamp, J.~Conway, R.~Guy, Winning Ways for Your Mathematical Plays,
  Academic Press, London, 1982, 697-734.

\bibitem{deBruijn*72}
N.~de~Bruijn, A {S}olitaire {G}ame and {I}ts {R}elation to a {F}inite {F}ield,
  Journal of Recreational Mathematics 5~(2) (1972) 133--137.

\bibitem{Avis-Deza-Onn*00}
D.~Avis, A.~Deza, S.~Onn, A combinatorial approach to the solitaire game,
  TIEICE: IEICE Transactions on Communications/Electronics/Information and
  Systems.

\bibitem{Lp-solve*06}
M.~Berkelaar, K.~Eikland, P.~Notebaert, lp\_solve version 5.5.0.6,
  \url{http://lpsolve.sourceforge.net/5.5/} (2006).

\bibitem{Deza-Onn*02}
A.~Deza, S.~Onn, Solitaire lattices, Graphs and Combinatorics 18~(2) (2002)
  227--243.

\bibitem{Avis-Deza*01}
D.~Avis, A.~Deza, On the binary solitaire cone, Discrete Applied Mathematics
  115~(1) (2001) 3--14.

\end{thebibliography}

\end{document}